\documentclass[a4paper,11pt]{article}
\usepackage[english]{babel}
\usepackage{a4wide}
\usepackage[utf8x]{inputenc}
\usepackage{hyperref}
\usepackage{setspace}
\usepackage{cite}
\usepackage{color}
\usepackage{amssymb}
\usepackage{amsmath}
\usepackage{amsthm}
\usepackage{graphicx}
\usepackage{enumerate}
\usepackage{dsfont}
\newtheorem{rmq}{Remark}
\newtheorem{thm}{Theorem}

\newcommand{\eps}{\varepsilon}
\newcommand{\R}{\mathbb{R}}

\newtheorem{lemme}{Lemma}

\title{Global existence for a damped wave equation and convergence towards a solution of the Navier-Stokes problem}
\author{Imène HACHICHA \footnote{Laboratoire Analyse et Probabilités, Université d'Evry-Val d'Essonne 23 bd de France, 91037 Evry Cedex. E-mail : imene.hachicha@univ-evry.fr}}
\date{}
\begin{document}

\newcommand{\tmax}{T^{\rm{max}}}
\newcommand{\Ed}{E^\delta}
\newcommand{\ddt}[1]{\frac{d #1}{dt}}
\newcommand{\ddx}[1]{\frac{d#1}{dx}}
\newcommand{\ddr}[1]{\frac{d#1}{dr}}
\newcommand{\be}{\begin{equation}}
\newcommand{\bea}{\begin{eqnarray}}
\newcommand{\beann}{\begin{eqnarray*}}
\newcommand{\ee}{\end{equation}}
\newcommand{\eea}{\end{eqnarray}}
\newcommand{\eeann}{\end{eqnarray*}}
\renewcommand{\P}{\mathbb{P}}
\maketitle

\abstract{In two and three space dimensions, and under suitable assumptions on the initial data, we show global existence for a damped wave equation which approaches, in some sense, the Navier-Stokes problem. The proofs are based on a refinement of the energy method in \cite{brenier}.

In this paper, we improve the results of \cite{brenier} and \cite{paicu}. We relax the regularity of the initial data of the former, even though we still use energy methods as a principal tool. Regarding \cite{paicu}, the improvement consists in the simplicity of the proofs and in requiring less regularity for the convergence to the Navier-Stokes problem. Indeed, the convergence result we obtain is near-optimal regularity.\\
\\
\\}

\section{Introduction}
First, let us recall the Navier-Stokes equations that describe the motion of a viscous, homegeneous and incompressible fluid
\be\left\{\begin{array}{ccc} \label{ns}
 \partial_t v - \Delta v = - \nabla : v \otimes v  - \nabla p  \\
 \nabla . v  =  0 \\
 v|_{t=0} = v_0 ,
\end{array}\right.\ee
where $v$ is the velocity of the fluid and $p$ the pressure, assumed to vanish at infinity, in the sense that
\be
\dfrac{1}{|B(0,R)|} \int_{B(0,R)}\!\!\!\!\!\! p ~~\longrightarrow ~0~,~~ R \rightarrow + \infty.
\ee
\\

The hyperbolic version of the Navier-Stokes equations studied here has been obtained after relaxation of the Euler equations and rescaling variables (see \cite{brenier} and references therein for the details):
\be\left\{\begin{array}{ccc} \label{nlw}
\eps \partial_{tt} u^\eps + \partial_t u^\eps - \Delta u^\eps = - \nabla : u^\eps \otimes u^\eps - \nabla p^\eps  \\
 \nabla . u^\eps  =  0 \\
(u^\eps, \partial_t u^\eps) |_{t=0} = (u_0^\eps , u_1^\eps) .
\end{array}\right.\ee
This type of approximation is the same as the hyperbolic perturbation of the heat equation introduced and investigated by Cattaneo in the fifties (see \cite{cattaneo} and further works).\\

\par
In \cite{brenier}, Brenier, Natalini and Puel proved global existence and uniqueness for the perturbed Navier-Stokes equation \eqref{nlw} with initial data in $H^2(\mathbb{T}^2 )^2 \times H^1(\mathbb{T}^2 )^2$, where $\mathbb{T}^2$ is the unit periodic square $\R^2 / \mathbb{Z}^2$, if $\eps < \eps_0$ and $\eps_0$ depends only on the initial data $(u_0^\eps , u_1^\eps )$.
Moreover, they proved the convergence of the solution to \eqref{nlw} towards a smooth solution to \eqref{ns} in the $L^\infty (]0,T[ ; L^2(\R^2))$ norm, provided that $v_0$ is smooth.\\

\par
In \cite{paicu}, Paicu and Raugel considered the same relaxed model \eqref{nlw}. In their paper, they proved, again if $\eps$ is small enough, global existence and uniqueness results for \eqref{nlw} with significantly improved regularity for the initial data. In fact, they need only $H^1(\R^2)^2 \times L^2(\R^2)^2$ regularity, thanks to a Strichartz estimate. In the three-dimensional case, they state a global existence result under a natural smallness condition on the initial data in $H^{1+\delta}(\R^3)^3 \times H^{\delta}(\R^3)^3$, for $\delta >0$ and $\eps$ small enough. Furthermore, they announced improved error estimates in the $H^1(\R^2)^2$ norm for $u^\eps -v$ and in the $L^2(\R^2)^2$ norm for $\partial_t \left( t (u^\eps -v) \right)$ in the two-dimensional case, where $u^\eps$ is as above and $v$ is the solution to \eqref{ns} with $v_0 = u_0^\eps \in H^1 (\R^2)^2$ or $v_0 = u_0^\eps \in H^{1+\delta} (\R^3)^3$.\\

\par
In two and three space dimensions and under suitable smallness assumptions on the initial data, we prove global existence for the Cauchy problem \eqref{nlw} in $\dot H^{\frac{n}{2}+\delta} \cap \dot H^{\frac{n}{2} -1 +\delta} (\R^n)^n$ in dimensions $n=2$ and $n=3$.

Moreover, for all positive $T$, we prove the convergence in the $L^\infty (]0,T[ ; \dot H^{\frac{n}{2}-1}(\R^n)^n)$ norm, with $n=2,3$, of solutions to \eqref{nlw} towards solutions to the Navier-Stokes problem \eqref{ns} with initial data $v_0 \in H^{\frac{n}{2}-1+s}(\R^n)^n$ and $s>0$.
More precisely, we prove the following two theorems.
\begin{thm} \label{th2d}
 Let $T\geq 0$ and $0< s,\delta <1$. Let $v_0 \in H^s (\R^2)^2$ be a divergence-free vector field and $(u_0^\eps , u_1^\eps) \in H^{1+\delta}(\R^2)^2 \times H^{\delta}(\R^2)^2$ be a sequence of initial data for problem \eqref{nlw}. Assume
\be\left\{\begin{array}{ccc} \label{hyp_th2d}
\|u_0^\eps -v_0\|_{L^2} + \eps \|u_1^\eps\|_{L^2} + \eps^{\frac{1}{2}} \|u_0^\eps\|_{\dot H^1} + \eps^{\frac{1+\delta}{2}} \|u_0^\eps \|_{\dot H^{1+\delta}} + \eps^{\frac{\delta}{2}} \|u_0^\eps\|_{\dot H^\delta}= \mathcal{O} \left(\eps^{\frac{s}{2}}\right) \\
 \eps^{1+\frac{\delta}{2}} \|u_1^\eps \|_{\dot H^\delta} = o \left(1\right).
\end{array}\right.\ee
Then, for $\eps$ small enough, there exists a global solution $u^\eps$ to system \eqref{nlw} that converges, when $\eps$ goes to $0$, in the  $L^\infty(]0,T[; L^2(\R^2)^2) $ norm, towards the unique solution $v$ to the incompressible Navier-Stokes equations \eqref{ns}, with $v_0$ as initial data. Moreover, there exists a constant $C_T$, depending only on $T$ and $v$,  such that
\be
\sup_{t\in [0,T]} \int_{\R^2} |u^\eps -v|^2 \, dx \leq C_T \eps^{\left( \frac{s}{2} \right)^-}.
\ee
\end{thm}
\begin{thm} \label{th3d}
Let $T\geq 0$ and $0< s,\delta <1$. Let $v_0 \in H^{s+\frac{1}{2}} (\R^3)^3$ be a divergence-free vector field and $(u_0^\eps , u_1^\eps) \in H^{\frac{3}{2}+\delta}(\R^3)^3 \times H^{\frac{1}{2}+\delta}(\R^3)^3$ be a sequence of initial data for problem \eqref{nlw} such that $\|u_0^\eps \|_{\dot H^{\frac{1}{2}}} < \frac{1}{16}$. Assume
\be\left\{\begin{array}{ccc} \label{hyp_th3d}
\|u_0^\eps -v_0\|_{ \dot H^{\frac{1}{2}} } + \eps \|u_1^\eps\|_{ \dot H^{\frac{1}{2}} } + \eps^{\frac{1}{2}} \|u_0^\eps\|_{\dot H^{\frac{3}{2}}}+  \eps^{\frac{1+\delta}{2}} \|u_0^\eps \|_{\dot H^{\frac{3}{2}+\delta}} + \eps^{\frac{\delta}{2}} \|u_0^\eps\|_{\dot H^{\frac{1}{2}+\delta}} = \mathcal{O} \left(\eps^{\frac{s}{2}}\right) \\
 \eps^{1+\frac{\delta}{2}} \|u_1^\eps \|_{\dot H^{\frac{1}{2}+\delta}} = o \left(1\right). 
\end{array}\right.\ee
Then, for $\eps$ small enough, there exists a global solution $u^\eps$ to system \eqref{nlw} that converges, when $\eps$ goes to $0$, in the  $L^\infty(]0,T[;  \dot H^{\frac{1}{2}}(\R^3)^3) $ norm, towards the unique solution $v$ to the incompressible Navier-Stokes equations \eqref{ns}, with $v_0$ as initial data. Moreover, there exists a constant $C_T$, depending only on $T$ and $v$, such that
\be
\sup_{t\in [0,T]}  \int_{\R^3} | \Lambda^{\frac{1}{2}} (u^\eps -v)|^2 \, dx \leq C_T \eps^{\left( \frac{s}{2} \right)^-}.
\ee
\end{thm}
\begin{rmq}
As a consequence of the assumptions $\|u_0^\eps \|_{\dot H^{\frac{1}{2}}} < \frac{1}{16}$ and $\|u_0^\eps -v_0\|_{ \dot H^{\frac{1}{2}} }  = \mathcal{O} \left(\eps^{\frac{s}{2}}\right) $, we obtain the smallness of $\|v_0 \|_{\dot H^{\frac{1}{2}}}$, which is a necessary condition to the existence of global solutions to the Navier-Stokes equation in $\R^3$.
\end{rmq}

\begin{rmq}
 We prove the convergence for initial data $v_0 \in H^{s+\frac{n}{2}-1}(\R^n)^n$, where $s>0$ and $H^{\frac{n}{2}-1}(\R^n)^n$ is the critical space for the Navier-Stokes equation.
\end{rmq}
\vspace*{1cm}

Before going further, we should check up that, given $v_0$, there exists a couple of functions $(u_0^\eps, u_1^\eps)$ satisfying all assumptions  \eqref{hyp_th2d} in Theorem \ref{th2d}, for example. First, since $u_1^\eps$ is not involved in the condition $\|u_0^\eps -v_0\|_{L^2}=\mathcal{O} \left(\eps^{\frac{s}{2}}\right) $, we can take $u_1^\eps \equiv 0$. Then let $u_0^\eps$ be defined by
$$
 \widehat{u_0^\eps}(\xi) = \widehat{v_0}(\xi)~ \mathds{1}_{|\xi| < \frac{1}{\sqrt{\eps}}}.
$$
Now a Bernstein inequality (see \cite{pglr} page 24) gives
$$
 \|u_0^\eps \|_{\dot H^\sigma} \leq \left( \sqrt{\eps} \right)^{s-\sigma} \|v_0\|_{\dot H^s}
$$
and easy calculations lead to the Jackson inequality
$$
 \|u_0^\eps -v_0\|_{L^2} \leq \left( \sqrt{\eps} \right)^s \|v_0 \|_{\dot H^s} .
$$
Therefore, the conditions \eqref{hyp_th2d} are fulfilled by at least this particular choice of $(u_0^\eps , u_1^\eps)$. The same arguments lead to the existence of a couple of functions satisfying the assumptions \eqref{hyp_th3d} in Theorem \ref{th3d}.
\vspace*{1cm}

The energy method we use in this paper is inspired by \cite{brenier} and refined in order to take into account the loss of regularity of the initial data. Therefore, we frequently have to use tame estimates and interpolations; this makes the proofs more technical than in \cite{brenier}. Moreover, we introduce a new energy, inspired from the classical one in \cite{brenier, paicu}, which is more convenient for the spaces we work in. Let us notice that we do not use any Strichartz estimate in this paper so the proofs are easier to understand than those in \cite{paicu}. Nevertheless, we lose the natural Strichartz regularity for the global existence results.\\
\vspace*{1cm}

This paper is organized as follows. In the next section, we treat the two-dimensional case and prove Theorem \ref{th2d}. In the last section, we adapt the results to the three-dimensional case and prove Theorem \ref{th3d}.
\section{The two-dimensional case : proof of Theorem \ref{th2d}} \label{dimension2}
This section includes two subsections. In the first one, we prove global existence for \eqref{nlw} using a fixed point method. Then, in subsection 2, we prove the convergence of this global solution to a solution to the Navier-Stokes equation, which is the last part of the statement of Theorem \ref{th2d}.\\

\textbf{Notation:} In the following, let $\mathbb{P}$ denote the Leray projection that maps a vector field to its zero-divergence part and, to alleviate the notations, let $L_T^p  \mathcal{E}$ denote the space $L^p \left( (0,T) ; \mathcal{E} \right) $.
\subsection{Global existence in $\dot H^{1+\delta}(\R^2) \times \dot H^\delta (\R^2)$}
Since $\mathbb{P}$ is a convolution operator on $\R^2$, it commutes with the differential operators. So, applying $\mathbb{P}$ to problem \eqref{nlw}, we obtain the damped (nonlinear) wave equation
\be
 (NLW_\eps) \left\lbrace
\begin{array}{rcl}
\eps \partial_{tt} u^\eps + \partial_t u^\eps - \Delta u^\eps & = & - \mathbb{P} \nabla :u^\eps \otimes u^\eps \\
(u^\eps , \partial_t u^\eps )|_{t=0} & = & (u_0^\eps , u_1^\eps ) .
\end{array}
\right.
\ee
which we consider with $(u_0^\eps , u_1^\eps ) \in H^{1+\delta}(\R^2) \times H^\delta (\R^2)$. \\
It is more convenient to study the same equation with parameter $\eps = 1$ :
$$
 (NLW) \left\lbrace
\begin{array}{rcl}
\partial_{tt} u + \partial_t u - \Delta u & = & - \mathbb{P} \nabla :u \otimes u \\
(u , \partial_t u )|_{t=0} & = & (u_0 , u_1 ) \in H^{1+\delta}(\R^2) \times H^\delta (\R^2) .
\end{array}
\right.
$$
In this purpose, let us set
\be 
\label{scaling}
u^\eps (\tau, y) = \dfrac{1}{\sqrt{\eps}} u(\dfrac{\tau}{\eps} , \dfrac{y}{\sqrt{\eps}} ).
\ee
This scaling transforms system $(NLW_\eps)$ into system $(NLW)$ with initial data
\beann
u_0 (x) & = & \sqrt{\eps} u_0^\eps \left( \sqrt{\eps} x \right) \\
u_1 (x) & = & \eps^{\frac{3}{2}} ~ u_1^\eps \left( \sqrt{\eps} x \right)
\eeann
As usual, let us denote by $\Box$ the D'Alembert operator: $\Box =  \partial_{tt} - \Delta $ and rewrite $(NLW)$ as follows.
\be
 (NLW) \left\lbrace
\begin{array}{rcl}
\Box u & = & -\partial_t u - \mathbb{P} \nabla : u \otimes u  =: F(u) \\
(u , \partial_t u )|_{t=0} & = & (u_0 , u_1 ) \in H^{1+\delta}(\R^2) \times H^\delta (\R^2).
\end{array}
\right.
\ee
Duhamel's formula for $u$ solution to $(NLW)$ is then
\begin{equation}
\label{duhamelondes}
u(t) = \cos (t \Lambda ) u_0 + \frac{\sin (t \Lambda )}{\Lambda} u_1 + \int_0^t \frac{\sin \left( (t-s) \Lambda \right)}{\Lambda} F\left( u(s) \right) \,ds =: \phi \left( u \right) (t)~,
\end{equation}
where $\Lambda$ is the differentiation operator $\Lambda = \sqrt{-\Delta}$, defined in Fourier variables by
$$
 \widehat{\Lambda f}(\xi) = |\xi| \hat f (\xi) .
$$
The proof of local existence in $\dot H^{1+\delta} \cap \dot H^\delta (\R^2, \R^2)$ is standard and, therefore, omitted.\\
We obtain a local solution to $(NLW)$, defined on $[0,T[$, for all positive $T$ such that
\begin{equation}
\label{temps_existence}
T  \leq \frac{1}{12 + 48 C \left( \|u_0 \|_{\dot H^{1+\delta}(\R^2)} + \|u_0\|_{\dot H^\delta (\R^2)} + \|u_1 \|_{\dot H^\delta(\R^2)} \right) } .
\end{equation}
In particular, as long as $\|u (t) \|_{X} := \|u (t) \|_{\dot H^{1+\delta}(\R^2)} + \|u (t) \|_{\dot H^{\delta}(\R^2)} + \| \partial_t u (t) \|_{\dot H^\delta(\R^2)} $ remains bounded, we can iterate the fixed point argument and extend the solution. 
\paragraph{Globalization} \label{globalisation2d}
Let us resume the initial equation $(NLW_\eps )$. One can easily check that
\be
\|u (t) \|_{X} = \eps^{\frac{\delta}{2}} \left( \sqrt{\eps} \|u^\eps \|_{\dot H^{1+\delta}(\R^2)} + \|u^\eps \|_{\dot H^{\delta}(\R^2)} + \eps \| \partial_t u^\eps  \|_{\dot H^\delta(\R^2)} \right) .
\ee
\\
\\
Then, for all non negative real $\delta$ and all non negative $t$, we define the energy
\be
 E^\delta_\eps (t) = \int_{\mathbb{R}^2} \left( \frac{1}{2}\vert \Lambda ^\delta (u^\eps + \varepsilon \partial_t u^\eps ) \vert ^2 + \frac{\varepsilon^2}{2} \vert \Lambda^\delta \partial_t u^\eps \vert^2 + \varepsilon \vert \Lambda^{1+\delta } u^\eps \vert^2 \right) \, d x
\ee
so that we have
\be
 \eps^{\frac{\delta}{2}} \left( \sqrt{\eps} \|u^\eps \|_{\dot H^{1+\delta}(\R^2)} + \|u^\eps \|_{\dot H^{\delta}(\R^2)} + \eps \| \partial_t u^\eps \|_{\dot H^\delta(\R^2)} \right)  \leq C \eps^{\frac{\delta}{2}} \sqrt{E^\delta_\eps}.
\ee
In this subsection, we shall prove the following inequality that yields the globality of the solution.
\begin{equation} \label{globalisation}
 E^\delta_\eps (t) \leq C_0 \eps^{-\delta}
\end{equation}
for all $t>0$, where $C_0$ is a universal constant. \\
\\
First, let us point out that $\dot{H}^{1+\delta} \cap L^\infty (\R^2)$ is an algebra and that the product estimate
\begin{equation}
\label{estimation_douce}
 \left\Vert fg \right\Vert_{\dot{H}^{1+\delta}(\R^2)} \leq C_1 \left( \left\Vert f \right\Vert_{\dot{H}^{1+\delta}} \left\Vert g \right\Vert_{L^\infty} + \left\Vert g \right\Vert_{\dot{H}^{1+\delta}} \left\Vert f \right\Vert_{L^\infty} \right).
\end{equation}
holds (see \cite{alinhac}) for all functions $f,g \in\dot{H}^{1+\delta} \cap L^\infty (\R^2)$. Moreover, we know that the homogeneous Besov\footnote{For definitions and properties of the Besov spaces, see the book by P.-G. Lemarié-Rieusset \cite{pglr}} space $\dot{B}_{2,1}^1(\mathbb{R}^2)$ embeds into $L^\infty (\mathbb{R}^2)$ and, interpolating, we obtain 
\begin{equation}
\label{interp_besov}
 \|f\|_{L^\infty} ~ \leq ~ \tilde C \|f\|_{\dot{B}_{2,1}^1} ~ \leq ~ C_2 \left\Vert f \right\Vert_{\dot{H}^\delta}^\delta . \left\Vert f \right\Vert_{\dot{H}^{1+\delta}}^{1-\delta}.
\end{equation}
\begin{lemme} \label{lemme_controle_energie}
Assume the following, when $\eps$ goes to zero :
$$
 (H) \left\lbrace \begin{array}{ll}
                   i)  &  \eps^{\frac{1+\delta}{2}} \|u_0^\eps\|_{\dot H^{1+\delta}} + \eps^{\frac{\delta}{2}} \|u_0^\eps\|_{\dot H^{\delta}} = o(1) \\
		   ii) &  \eps^{\frac{1}{2}} \|u_0^\eps\|_{\dot H^{1}} + \eps  \|u_1^\eps\|_{L^2} =o(1) .
                  \end{array}
\right.
$$
Let us define $0 \leq T \leq \tmax_\eps$ by
\be
 T = \sup \left\lbrace 0 \leq \tau \leq \tmax_\eps ~ :~ \forall ~ t \in [0, \tau[ , ~\|u^\eps  (t)\|_{L^\infty} < \dfrac{1}{8~ C_1 \sqrt{\eps}} \right\rbrace .
\ee
Then, for $\eps$ small enough, there exists a large number $N$, only depending on $\delta$ and $\|u_0^\eps\|_{L^2}$ (which is arbitrary), such that, for all $0 \leq t \leq T$,
\begin{equation} \label{controle_energie}
 E^\delta_\eps (t) \leq E^\delta_\eps (0) \left(2 \|u_0^\eps\|_{L^2}^2 +1\right)^N .
\end{equation}
\end{lemme}
\proof
Let us compute the (time) derivative of $E^\delta_\eps$. 
\beann
 \ddt{ E^\delta_\eps} (t) & = & \int \Lambda^\delta (\varepsilon \partial_{tt} u^\eps + \partial_t u^\eps ).\Lambda^\delta (u^\eps + \varepsilon \partial_t u^\eps) + \varepsilon^2 \Lambda^\delta (\partial_t u^\eps).\Lambda^\delta (\partial_{tt} u^\eps) + 2 \varepsilon \Lambda^{\delta+1} (\partial_t u^\eps).\Lambda^{\delta+1} (u^\eps) \\
& = & \int \Lambda^\delta (\varepsilon \partial_{tt} u^\eps + \partial_t u^\eps - \Delta u^\eps) .\Lambda^\delta (u^\eps + 2 \varepsilon \partial_t u^\eps) - \varepsilon \vert\Lambda^\delta \partial_t u^\eps \vert^2 -  \vert \Lambda^{\delta +1} u^\eps \vert^2 \\
& = & \int - \Lambda^\delta (u^\eps. \nabla u^\eps).\Lambda^\delta u^\eps - 2 \varepsilon \Lambda^\delta(u^\eps . \nabla u^\eps).\Lambda^\delta \partial_t u^\eps - \varepsilon \vert \Lambda^\delta \partial_t u^\eps \vert^2 - \vert \Lambda^{\delta+1} u^\eps \vert^2 
\eeann
Performing a classical Young inequality $2ab \leq a^2 + b^2$ and rearranging the terms, we obtain
\beann
\ddt{ E^\delta_\eps} (t) 
& \leq & - \int  \Lambda^\delta (u^\eps . \nabla u^\eps).\Lambda^\delta u^\eps + \left(\varepsilon \int \vert \Lambda^\delta(u^\eps . \nabla u^\eps) \vert^2 - \int \vert \Lambda^{\delta+1} u^\eps \vert^2 \right) \\
\eeann
Now, let us recall that
\be
  u^\eps.\nabla u^\eps  = \sum_{i=1}^2 u^\eps_i \partial_i u^\eps  = \sum_{i=1}^2 \partial_i (u^\eps_i u^\eps) 
\ee
since $u^\eps$ is divergence-free. Now, \eqref{estimation_douce} yields
\be
 \| u^\eps.\nabla u^\eps \|_{\dot H^\delta} \leq  \sum_{i=1}^2 \|\partial_i (u^\eps_i u^\eps) \|_{\dot H^\delta} \leq 4 C_1 \|u^\eps\|_{L^\infty} \|u^\eps \|_{\dot H^{1+\delta}}
\ee
Thus
\begin{eqnarray*}
 \left| \int  \Lambda^\delta (u^\eps . \nabla u^\eps).\Lambda^\delta u^\eps \,dx \right| & \leq & \lVert u^\eps . \nabla u^\eps \rVert_{\dot{H}^\delta}  \lVert u^\eps \rVert_{\dot{H}^\delta}\\
& \leq & 4 C_1 \|u^\eps\|_{L^\infty} \|u^\eps \|_{\dot H^{1+\delta}} \lVert u^\eps \rVert_{\dot{H}^\delta}\\
& \leq & 4 C_1 C_2\lVert u^\eps  \rVert_{\dot{H}^\delta}^\delta \lVert u^\eps  \rVert_{\dot{H}^\delta} \lVert u^\eps  \rVert_{\dot{H}^{1+\delta}}^{2-\delta} .
\end{eqnarray*}
The interpolation inequality
\be
 \|u^\eps \|_{\dot{H}^\delta} \leq C_3 \|u^\eps \|_2^{1-\delta} \lVert u^\eps \rVert_{\dot{H}^1}^\delta 
\ee
yields finally
\be
 \left| \int  \Lambda^\delta (u^\eps . \nabla u^\eps).\Lambda^\delta u^\eps \,dx \right| \leq 4 C_1 C_2 C_3 \| u^\eps \|_2^{1-\delta} \left( \lVert u^\eps  \rVert_{\dot{H}^1} \lVert u^\eps  \rVert_{\dot{H}^\delta} \right)^\delta \lVert u^\eps  \rVert_{\dot{H}^{1+\delta}}^{2-\delta}.
\ee
Besides, since we assume $i)$ and using inequality \eqref{interp_besov}, we can write, for $\eps$ small enough,
\begin{equation}
\label{controle_norme_infinie}
 \| u_0^\eps \|_{L^\infty} \leq \dfrac{1}{8 C_1 \sqrt{\varepsilon}}  .
\end{equation}
Now, by continuity of the (local) solution $u^\eps$ with respect to $t$, we deduce that $T>0$ and the inequality 
\begin{eqnarray}
\eps \int \vert \Lambda^\delta(u^\eps . \nabla u^\eps ) \vert^2 - \int \vert \Lambda^{\delta+1} u^\eps  \vert^2 & \leq & \left( 16 C_1^2 \varepsilon \| u^\eps  \|_{L^\infty}^2 -1 \right) \int |\Lambda^{\delta+1} u^\eps  |^2 \nonumber \\
& \leq & -\dfrac{1}{4} \| u^\eps \|_{\dot{H}^{1+\delta}}^2  \label{est}
\end{eqnarray}
holds on $[0,T[$. \\
Hence, for all $0 \leq t <T$,  
\be
 \frac{d}{dt} E^\delta_\eps (t) ~ \leq ~ -\dfrac{1}{4} \|u^\eps \|_{\dot{H}^{1+\delta}}^2 + 4 C_1 C_2 C_3  \|u^\eps \|_2^{1-\delta} \left( \lVert u^\eps  \rVert_{\dot{H}^1} \lVert u^\eps  \rVert_{\dot{H}^\delta} \right)^\delta \lVert u^\eps  \rVert_{\dot{H}^{1+\delta}}^{2-\delta}.
\ee
Now, consider the second term on the right hand side and use Young inequality
\be
 ab \leq \frac{\delta}{2}~ a^{\frac{2}{\delta}} + \frac{2-\delta}{2}~b^{\frac{2}{2-\delta}} 
\ee
with $b=2^{\delta -2} \|u^\eps \|_{\dot{H}^{1+\delta}}^{2-\delta}$. We obtain
\beann
\frac{d}{dt} E^\delta_\eps (t) & \leq & \dfrac{1}{4} \left( \frac{2-\delta}{2} -1 \right) \| u^\eps \|_{\dot{H}^{1+\delta}}^2 +
\frac{\delta}{2} \left( 2^{2-\delta} 4 C_1 C_2 C_3 \| u^\eps \|_2^{1-\delta} \| u^\eps \|_{\dot{H}^1}^\delta \| u^\eps \|_{\dot{H}^\delta}^\delta \right)^{\frac{2}{\delta}} \\
& \leq & \frac{\delta}{2}~ 2^{2\frac{2-\delta}{\delta}} (4 C_1 C_2 C_3 )^{\frac{2}{\delta}} \| u^\eps \|_2^{2\frac{1-\delta}{\delta}} \| u^\eps \|_{\dot{H}^1}^2 \| u^\eps \|_{\dot{H}^\delta}^2 .
\eeann
Then the inequality
\be
 \| a+b \|^2 + \| a-b \|^2  = 2 \|a\|^2 + 2 \|b\|^2 \geq 2 \| b \|^2
\ee
with $a=\dfrac{u^\eps }{2} + \eps \partial_t u^\eps $ and $b=\dfrac{u^\eps }{2}$ yields, for all non negative $\delta$, the estimate 
\begin{equation} 
\label{controle_Hdelta}
E^\delta_\eps (t) ~ \geq ~ \| u^\eps  + \eps \partial_t u^\eps  \|_{\dot{H}^\delta}^2 + \| \eps \partial_t u^\eps  \|_{\dot{H}^\delta}^2 ~ \geq ~ \frac{1}{2} \| u^\eps \|_{\dot{H}^\delta}^2 
\end{equation}
from which we deduce
\begin{equation}
\frac{d}{dt} E^\delta_\eps (t) \leq \delta ~ 2^{2\frac{2-\delta}{\delta}} (4 C_1 C_2 C_3 )^{\frac{2}{\delta}} \| u^\eps \|_2^{2\frac{1-\delta}{\delta}} \| u^\eps \|_{\dot{H}^1}^2 E^\delta_\eps (t).
\end{equation}
\\

Instead of showing that $t \mapsto E^\delta_\eps (t)$ is decreasing, we will prove that $E^ \delta (E^0_\eps +1)^N$ decreases if $N$ is large enough. Even though this new energy does not seem natural, it is more convenient than $E^ \delta +N E^0_\eps $ since the latter would have required a smallness assumption on the $L^2$ norm of $u_0^\eps$. But this would yield a weaker result than one might hope as the Navier-Stokes equation with initial data in $L^2(\R^2)$ unconditionally has a global solution.\\

First of all, let us notice that $\int_{\R^n} f . (f. \nabla f ) =0$ for all divergence-free function $f$ and that (see inequality \eqref{est})
\be
\ddt{E^0_\eps}(t) ~ \leq ~ -\frac{1}{4} \| u^\eps (t) \|_{\dot{H}^1}^2 .
\ee
We use it in \eqref{controle_derivee_energie} below.
\bea
\frac{d}{dt} \left[ (E^0_\eps +1)^N E^\delta_\eps  \right] & = & N \frac{d E^0_\eps}{dt} (E^0_\eps+1)^{N-1} E^\delta_\eps + (E^0_\eps +1)^N \frac{d E^\delta_\eps }{dt} \nonumber \\
& \leq & \frac{1}{4} ~(E^0_\eps +1)^{N-1} E^\delta_\eps \| u^\eps \|_{\dot{H}^1}^2  \times \nonumber \\
& &  \times \left[ -N + \delta ~ (16 C_1 C_2 C_3 )^{\frac{2}{\delta}} (E^0_\eps +1) \| u^\eps \|_2^{2 \frac{1-\delta}{\delta}}  \right]. \label{controle_derivee_energie}
\eea
Under the assumption $iii)$, we have
\beann
E^0_\eps(0) & = & \int_{\R^2} \left( \dfrac{|u_0^\eps  + \eps u_1^\eps  |^2 + \eps^2 |u_1^\eps |^2}{2} + \eps |\nabla u_0^\eps  |^2 \right) dx \\
& \leq & \|u_0^\eps \|_2^2 + \dfrac{3\eps ^2}{2} \|u_1^\eps \|_2^2 + \eps \|u_0^\eps \|_{\dot{H}^1}^2 \\
& < & 2~ \|u_0^\eps \|_2^2 
\eeann 
if $\eps$ is small enough. Now, \eqref{controle_Hdelta} with $\delta = 0$ yields
\begin{equation}
\label{conditions_u0}
 \| u^\eps  (t)\|_2^2  ~ \leq ~ 2 E^0_\eps(t) ~ < ~ 4 \| u_0^\eps \|_2^2
\end{equation}
for all $t \in (0,T)$. We have obtained
\be
 \frac{d}{dt} \left[ (E^0_\eps +1)^N E^\delta_\eps  \right] ~ \leq ~ \dfrac{1}{4} (E^0_\eps +1)^{N-1} E^\delta_\eps \|u^\eps \|_{\dot{H}^1}^2  \left[ -N + \delta ~ (16 C_1 C_2 C_3)^{\frac{2}{\delta}} (4 \|u_0^\eps \|_2^2 )^{\frac{1-\delta}{\delta}}  (2 \|u_0^\eps \|_2^2 +1)  \right].
\ee
\\
We deduce that, for $N$ large enough (depending only on $\delta$ and $\|u_0^\eps \|_2$), the right hand side is negative and
\be
 E^\delta_\eps (t) \leq E^\delta_\eps (0) (E^0_\eps(0) +1)^N < E^\delta_\eps(0) (2 \|u_0^\eps \|_2^2 +1)^{N}
\ee
for all $0\leq t < T$. \qed \\
\\
\\
\par
We want to reiterate the reasoning, $i.e.$ to keep the control of $\|u^\eps (t)\|_{L^\infty}$. The aim of Lemma \ref{lemme2} is to ensure this control throughout the time.
\begin{lemme} \label{lemme2}
Assume the limit 
\be \label{hyp}
 \eps^{1+\frac{\delta}{2}} \|u^\eps_1\|_{\dot H^\delta} \longrightarrow 0~, ~~\eps \rightarrow 0
\ee
in addition to the assumptions $(H)$ in Lemma \ref{lemme_controle_energie}. Then the inequality
\begin{equation}
\label{controle_initial}
 E^\delta_\eps (0) ~ < ~ \frac{2^{-\delta}}{(16 C_1 C_2)^2} ~ \eps^{- \delta}  (2 \| u_0^\eps \|_2^2 +1)^{-N}.
\end{equation}
holds for $\eps$ small enough. Moreover, for all $t \in [0,T[$,
\begin{equation}  \label{ongardelecontrole}
\|u^\eps \|_{L^\infty} ~\leq ~ \frac{1}{16 C_1 \sqrt{\eps}} .
\end{equation}
\end{lemme}
\proof
Easy calculations show that \eqref{controle_initial} is true as a result of $(H)$ and \eqref{hyp}.
Now, let $t \in [0,T[$ and recall interpolation inequality \eqref{interp_besov} 
\be
 \|u^\eps\|_{L^\infty} \leq C_2 \|u^\eps\|_{\dot H^\delta}^\delta \|u^\eps\|_{\dot H^{1+\delta}}^{1-\delta}.
\ee
Thanks to inequality \eqref{controle_Hdelta}, we have 
\be \| u^\eps (t) \|_{L^\infty}  \leq  C_2 \sqrt{2 E^\delta_\eps (t)}^\delta \sqrt{\eps^{-1}E^\delta_\eps (t)}^{1-\delta} =  C_2 2^{\frac{\delta}{2}} \eps^{\frac{\delta - 1}{2}} \sqrt{E^\delta_\eps (t)} . \ee
Using \eqref{controle_energie} (conclusion of Lemma \ref{lemme_controle_energie}) and \eqref{controle_initial}, we obtain
\be \| u^\eps (t) \|_{L^\infty}  \leq  C_2 2^{\frac{\delta}{2}} \eps^{\frac{\delta - 1}{2}} \frac{2^{-\frac{\delta}{2}}}{16 C_1 C_2}\eps^{-\frac{\delta}{2}}  =  \frac{1}{16C_1} \eps^{-\frac{1}{2}} 
\ee
for all $0<t<T$. \qed \\
\\
\begin{rmq}
 Under the assumptions \eqref{hyp_th2d} in Theorem \ref{th2d}, the conditions $(H)$ in Lemma \ref{lemme_controle_energie} and \eqref{hyp} in Lemma \ref{lemme2} are fulfilled.
\end{rmq}
\par
Now, we shall prove that these estimations (on $\|u^\eps (t) \|_{L^\infty}$ and, consequently, on $E^\delta_\eps (t)$) remain true on the whole existence interval $[0,\tmax_\eps)$, where $\tmax_\eps$ is the existence time \eqref{temps_existence} given by the Picard iteration. 
We have already proved that $T > 0$. Assume $T < \tmax_\eps$. Then
\begin{equation}
\label{contradiction}
 \| u^\eps  (T)\|_{L^\infty} = \dfrac{1}{8 C_1 \sqrt{\eps}}.
\end{equation}
On the other hand, \eqref{ongardelecontrole} in Lemma \ref{lemme2} yields 
\be
 \|u^\eps  (T)\|_{L^\infty} \leq \dfrac{1}{16 C_1 \sqrt{\eps}} < \dfrac{1}{8 C_1 \sqrt{\eps}},
\ee
which contradicts \eqref{contradiction} so $T \geq \tmax_\eps$.
We deduce then from Lemma \ref{lemme_controle_energie} and Lemma \ref{lemme2} that $E^\delta_\eps$ satisfies inequality \eqref{globalisation} on the existence interval $[0,\tmax_\eps)$. Therefore the $(NLW_\eps)$ equation has a global solution.

\subsection{Convergence towards a solution to Navier-Stokes problem} \label{convergence2d}
\par
Let us recall that Brenier, Natalini and Puel showed in \cite{brenier} that, under suitable assumptions, the solutions to $(NLW_\eps)$, with initial data $(u_0^\eps , u_1^\eps) \in H^2 \times H^1 ( \mathbb{T}^2)$, converge to the solutions to $(NS)$ with smooth initial data $v_0$ when $\eps$ goes to $0$. The authors used the modulated energy method (see \cite{brenier} and references therein) to show an error estimate in the $L^\infty_T L^2(\mathbb{T}^2)$ norm, for all positive $T$.\\

In this section, we show a similar result with less regularity on the initial data $(u_0^\eps ,u_1^\eps)$ and $v_0$, in $\R^2$ instead of $\mathbb{T}^2$. Indeed, we prove that, under suitable assumptions (less restrictive than those in \cite{brenier}) and for any positive $\delta$, the solutions to $(NLW_\eps)$, with initial data $(u_0^\eps , u_1^\eps) \in H^{1+\delta} \times H^\delta ( \mathbb{R}^2)$, converge in the $L^2(\R^2)^2$ norm, when $\eps$ goes to $0$, to the solutions to $(NS)$ with initial data $v_0 \in \dot H^{s}(\R^2)^2$, $s>0$. We recall here that the critical space for the Navier-Stokes equations in $\R^2$ is $L^2$.\\
Because of the loss of regularity, the proof is not straightforward anymore. Interpolation and product estimates will help to get around the difficulties.
\\
\par
We shall start the proof like in \cite{brenier}, introducing the so-called Dafermos modulated energy, which is the total energy of the wave equation $(NLW_\eps)$, modulated by a divergence-free vector $(t,x)\mapsto v(t,x)$
\begin{equation} \label{dafermos}
 E_{\eps, v} (t) = \int_{\R^2} \left( \frac{1}{2} |u^\eps -v(t,x) +\eps \partial_t u^\eps |^2 + \frac{\eps^2}{2} |\partial_t u^\eps |^2 + \eps |\nabla u^\eps |^2 \right) \,dx.
\end{equation}
This energy satisfies the inequality
\begin{equation} \label{u-v}
 \int_{\R^2} |u^\eps -v |^2 \,dx \leq 4 E_{\eps, v}(t).
\end{equation}
Via a Gronwall inequality, we shall show that, for all positive $t$ such that $t<T$, the modulated energy $E_{\eps,v}(t)$ converges to $0$ when $\eps$ goes to $0$. In this purpose, let us compute the time derivative of $E_{\eps,v}$.\\

\begin{lemme}[see \cite{brenier}]
If $u^\eps$ and $v$ are divergence-free functions, the Dafermos modulated energy satisfies the identity
\begin{eqnarray}
 \frac{d}{dt} E_{\eps,v}(t) & = & \int v. \nabla : (u^\eps -v )\otimes (u^\eps -v) -\eps \int |\partial_t u^\eps + \nabla : (u^\eps \otimes u^\eps)|^2  \nonumber \\
& & ~~~- \eps \int \partial_t v . \partial_t u^\eps + \int \left( \eps |\nabla : (u^\eps \otimes u^\eps)|^2 - |\nabla (u^\eps -v)|^2 \right) \nonumber \\
& & ~~~  + \int (\partial_t v + v.\nabla v -\Delta v)(v-u^\eps). \label{soisnul}
\end{eqnarray}
\end{lemme}
\begin{rmq}
 This lemma is proved in \cite{brenier}.
\end{rmq}
From now on, let $v$ be a solution to the Cauchy problem
\begin{equation}
\label{nschaleur}
\partial _t v = \Delta v - \P \nabla : v\otimes v ~,~~ \nabla . v =0 ~,~~ v|_{t=0} = v_0 \in H^s (\mathbb{R}^2) ,
\end{equation}
where $s$ is a positive real, so that the last term of \eqref{soisnul} vanishes. Now, since the second term in the identity \eqref{soisnul} is negative, \eqref{soisnul} writes
\begin{eqnarray}
 \frac{d}{dt} E_{\eps,v}(t) & \leq & - \eps \int_{\R^2} \partial_t v . \partial_t u^\eps \, dx +\int_{\R^2} v. \nabla : (u^\eps -v )\otimes (u^\eps -v) \, dx  \nonumber \\
& & + \int_{\R^2} \left( \eps |\nabla : (u^\eps \otimes u^\eps)|^2 - |\nabla (u^\eps -v)|^2 \right) \, dx. \label{derivee_energie}
\end{eqnarray}
We shall treat the three terms on the right hand side in the following three subsections.\\
\\
First, let us recall the Duhamel formula for Navier-Stokes equation \eqref{nschaleur}:
\begin{equation}
 \label{duhamel}
v(t,.)=e^{t\Delta}v_0 - \int_0^t e^{(t-s)\Delta} \P \nabla : (v\otimes v)(s,.)\,ds .
\end{equation}
Computing the time derivative of \eqref{duhamel}, we obtain
\begin{equation}
 \label{dtv}
\partial _t v(t,.)=\Delta e^{t\Delta}v_0 - \P \nabla : (v\otimes v)(t,.) - \int_0^t \Delta e^{(t-s)\Delta} \P \nabla : (v\otimes v)(s,.)\,ds .
\end{equation}
\begin{rmq}
\label{v}
If $v_0 \in H^s (\R^2)$, then, for all positive $T$, the solution $v$ to \eqref{nschaleur} satisfies
\begin{equation}
 v \in L^2_T  H^{s+1}(\mathbb{R}^2) \cap L_T ^\infty  H^s (\mathbb{R}^2).
\end{equation}
\end{rmq}
\subsubsection{Estimating $\varepsilon \int_0^T\!\!\! \int_{\mathbb{R}^2} \partial_t u^\varepsilon \partial_t v \,dt\,dx $} \label{partie1convergence2d}
In order to bound the integral by a positive power of $\eps$, we shall find spaces $L^p_T  \dot H^\sigma (\R^2)$ and $L^q_T  \dot H^{\sigma'} (\R^2)$  containing respectively $\partial_t v$ and $\eps \partial_t u^\eps$. 
\paragraph{Step 1: Regularity of $\partial _t v$.}
\begin{lemme} \label{regdtv}
 Let $\theta \in [0,1[$ and $v$ be a solution to equation \eqref{nschaleur}. Then, for all $0< \tilde \eps < s$, 
\begin{equation} \label{interpolationdtv}
 \partial _t v \in L^p_T  \dot H^\sigma (\mathbb{R}^2) ,
\end{equation}
where $\dfrac{1}{p}=\dfrac{1-\theta}{2}+\theta = \dfrac{1+\theta}{2} $ and $\sigma = (1-\theta ) (s-1) + \theta (s-\tilde \varepsilon) = s-1 +\theta (1-\tilde \varepsilon)$.
\end{lemme}
\proof
Recall that 
\be
\partial _t v(t,.)=\Delta e^{t\Delta}v_0 - \nabla : (v\otimes v)(t,.) - \int_0^t \Delta e^{(t-s)\Delta} \nabla : (v\otimes v)(s,.)\,ds .
\ee
First, standard estimates yield
\be
 \Delta e^{t\Delta}v_0 \in L^2_T \dot H^{s+1} \cap L^1_T \dot H^{s-\tilde \eps} .
\ee
Let us consider now the term $\nabla : (v\otimes v)$. Using remark \ref{v} and a classical Sobolev embedding, we have
$v \in L^2_T  L^\infty \cap L^\infty_T  H^s$. Since $L^\infty \cap \dot H^s$ is an algebra (see \cite{alinhac} page 98), we deduce that
\begin{equation}
 (v \otimes v )_{i,j} ~ = ~ v_i v_j \in L^2_T  \dot H^s
\end{equation}
Consequently,
\begin{equation}
 \nabla : (v\otimes v) \in L^2_T  \dot H^{s-1}.
\end{equation}
Moreover, using that $v \in L^2_T  H^{s+1}(\mathbb{R}^2)$ and $H^{s+1}(\mathbb{R}^2)$ being an algebra (see \cite{alinhac}), we obtain
$ (v \otimes v )_{i,j} ~ = ~ v_i v_j \in L^1_T  H^{s+1}(\mathbb{R}^2)$, so that 
\[ \nabla : (v\otimes v) \in L^1_T  H^{s}(\mathbb{R}^2) . \]
Interpolating between $L_T^2 \dot H^{s-1}(\R^2)$ and $L^1_T \dot H^{s-\tilde \eps}$ and taking $\theta \in [0,1]$, we have 
\begin{equation}
 \Delta e^{t \Delta} v_0 ~,~~ \P \nabla : (v\otimes v) \in L^{\frac{2}{1+\theta}} \dot H^{s-1+\theta(1-\tilde \eps)}.
\end{equation}
The following theorem (proved in \cite{pglr} page 64), applied to $\P \nabla :(v\otimes v)$, allows us to conclude immediately that the integral term is also in $ L^{\frac{2}{1+\theta}} \dot H^{s-1+\theta(1-\tilde \eps)}$ for all $\theta \in [0,1)$.
\begin{thm}[Maximal $L^p  L^q$ regularity for the heat kernel]
\label{regmax}
 Let $A$ be an operator defined by \be f(t, x) \mapsto Af(t, x) = \int_0^t e^{(t-s) \Delta} \Delta f (s,x) \, ds .\ee
Then $A$ is bounded from $L^p_T  L^q(\mathbb{R}^d)$ to itself, for all reals $T >0$, $1 < p < \infty$ and $1 < q < \infty$.
\end{thm}
\qed
\paragraph{Step 2: Regularity of $\varepsilon \partial _t u^\varepsilon$.}
\begin{lemme} Let $0 \leq \theta \leq 1$ and $u^\eps$ be a solution to $(NLW_\eps)$. Then 
\begin{equation} \label{regdtu}
 \eps \partial_t u^\eps \in L^q_T  \dot H^{\sigma'}( \mathbb{R}^2 ) 
\end{equation}
where $\dfrac{1}{q}=\dfrac{1-\theta}{2}$ and $\sigma' = \theta \delta$ and \be \| \eps \partial_t u^\eps \|_{L^q_T  \dot H^{\sigma'}} = \mathcal{O} \left(\eps^{\frac{1-\theta}{2}-\frac{\theta \delta}{2}} \right).\ee
\end{lemme}
\proof
First, let us show that
\begin{equation} \label{regdtu1}
 \varepsilon \partial_t u^\varepsilon \in L^2_T  L^2( \mathbb{R}^2 , \varepsilon^{-1/2}dx) \cap L^\infty_T  \dot H^\delta( \mathbb{R}^2 , \varepsilon^{\delta/2}dx).
\end{equation}
In this purpose, let us consider the energy
\be E^0_\eps (t)=:\int_{\mathbb{R}^2} \left( \dfrac{1}{2} |u^\eps + \eps \partial _t u^\eps |^2 + \dfrac{\eps ^2}{2} |\partial _t u^\eps |^2 + \eps |\nabla u^\eps |^2 \right) \,dx \ee and differentiate it. A simple calculation gives us
\be
 \frac{d}{dt}E^0_\eps (t) + \eps \int_{\mathbb{R}^2}|\partial _t u^\eps + \nabla : (u^\eps \otimes u^\eps )|^2 \,dx + \int_{\mathbb{R}^2} \left( |\nabla u^\eps |^2 - \eps |\nabla :(u^\eps \otimes u^\eps)|^2 \right)\,dx =0.
\ee
Thanks to the control of the norm $\|u^\eps \|_{L^\infty}$ by $\dfrac{1}{2\sqrt{\eps}}$, the last term in the left hand side is positive and we obtain
\be
 E^0_\eps (T) + \int_0^T \!\!\! \int_{\mathbb{R}^2}\eps |\partial _t u^\eps + \nabla : (u^\eps \otimes u^\eps )|^2 \,dx\,dt + \int_0^T \!\!\!\int_{\mathbb{R}^2} \frac{1}{2} |\nabla u^\eps |^2 \,dx\,dt \leq E^0_\eps (0) \leq C_0 .
\ee
Therefore, we have that $\int_0^T \!\!\! \int_{\mathbb{R}^2}\eps |\partial _t u^\eps + \nabla : (u^\eps \otimes u^\eps )|^2 \,dx\,dt \leq C$ and $\int_0^T \!\!\!\int_{\mathbb{R}^2} \frac{1}{2} |\nabla u^\eps |^2 \,dx\,dt \leq C$ which gives
\be
 \eps \int_0^T \!\!\!\int_{\mathbb{R}^2} | \nabla : (u^\eps \otimes u^\eps )|^2 \,dx\,dt \leq \eps \|u^\eps\|_{L^\infty} \int_0^T \!\!\!\int_{\mathbb{R}^2}  |\nabla u^\eps |^2 \,dx\,dt \leq C .
\ee
Consequently,
\begin{eqnarray}
\| \sqrt{\eps} \partial_t u^\eps \|_{L^2_T  L^2}^2 & \leq & 2 \int_0^T \!\!\! \int_{\mathbb{R}^2}\eps |\partial _t u^\eps + \nabla : (u^\eps \otimes u^\eps )|^2 \,dx\,dt + 2 \eps \int_0^T \!\!\!\int_{\mathbb{R}^2} | u^\eps . \nabla  u^\eps |^2 \,dx\,dt \nonumber \\
& \leq & 2C + 2C = 4C,
\end{eqnarray}
\emph{i.e.} $\sqrt{\eps} \partial_t u^\eps \in L^2_T  L^2$ uniformly in parameter $\eps$. \\
\\
Moreover, directly from the expression of $E^\delta_\eps$ and from \eqref{globalisation}, we deduce that 
\be \|\partial_t u^\varepsilon\|_{L^\infty_T  \dot H^\delta( \mathbb{R}^2)} = \mathcal{O} \left( \varepsilon^{-1-\frac{\delta}{2}} \right).\ee
\\
Now, interpolating between the two spaces in \eqref{regdtu1}, we obtain \eqref{regdtu}. \qed
\\
\\
\\
We will now choose un appropriate $0< \theta <1$, such that 
\be
\varepsilon \int_0^T\!\!\! \int_{\mathbb{R}^2} \partial_t u^\varepsilon \partial_t v \,dt\,dx \longrightarrow 0 ~,~~ \eps \rightarrow 0. 
\ee
So we should have
\begin{equation} \label{theta}
 -1+\theta (1+\delta) <0 \Longleftrightarrow 0< \theta < \dfrac{1}{1+\delta},
\end{equation}
so that the power of $\eps$ controlling the integral is negative. Moreover, $L^q_T  H^{\sigma'}$ is the dual space of $L^p_T  H^\sigma$ if $\sigma' = - \sigma$ , $i.e.$
\be \theta = \frac{1-s}{1+\delta - \tilde \varepsilon}. \ee
This value of $\theta$ satisfies condition \eqref{theta} and gives us
\be
 \varepsilon \int_0^T \!\!\! \int_{\mathbb{R}^2} \partial_t u^\varepsilon \partial_t v \,dt\,dx = \mathcal{O} \left(\sqrt{\varepsilon} ^{(s(1+\delta ) - \tilde \varepsilon )/(1+\delta - \tilde \varepsilon)} \right) = \mathcal{O} \left(\varepsilon ^{\nu} \right),
\ee
with $\nu = \dfrac{s(1+\delta ) - \tilde \varepsilon }{2(1+\delta - \tilde \varepsilon)} >0$. Notice that $0<\nu < \frac{s}{2}$ and that $\nu$ can be chosen arbitrarily close to $\frac{s}{2}$.
\subsubsection{Estimating $ \int_0^T\!\!\! \int_{\mathbb{R}^2} v. \nabla : (u^\eps - v ) \otimes (u^\eps -v) \,dt\,dx $}
\begin{lemme}
 Let $v$ and $u^\eps$ be as above. We have
\begin{equation}
 \label{est1}
\int_0^T\!\!\! \int_{\mathbb{R}^2} v. \nabla : (u^\eps - v ) \otimes (u^\eps -v) \,dt\,dx \leq 2 \int_0^T \|v\|_{BMO}^2 E_{\eps, v} (t) \,dt + \frac{1}{2} \int_0^T \|\nabla (u^\eps -v)\|_{L^2}^2 .
\end{equation}
\end{lemme}
\proof
Let us recall that $v(t) \in \dot H^1 (\mathbb{R}^2) \subset BMO (\R^2)$ (see remark \ref{v}) and that the solution of the wave equation $(NLW_\eps)$, $u^\eps (t)$, is an $L^2 (\mathbb{R}^2)$ divergence-free vector. The following theorem applies.
\begin{thm}[div-curl lemma, see \cite{clms}]
\label{divcurl}
 Let $f$ be an $L^2 (\mathbb{R}^2)$ divergence-free vector and $g\in \dot H^1 (\mathbb{R}^2)$. Then \be f.\nabla g \in \mathcal{H}^1 (\mathbb{R}^2),\ee where $\mathcal{H}^1 (\mathbb{R}^2)$ is the Hardy space constructed on $L^1 (\mathbb{R}^2)$.
\end{thm}
Using the duality of Hardy space $\mathcal{H}^1 (\mathbb{R}^2)$ and $BMO (\mathbb{R}^2)$ (proved in \cite{stein}), we write
$$\int_{\mathbb{R}^2} v. \nabla : (u^\eps - v ) \otimes (u^\eps -v) \,dx  \leq  c' \|v\|_{BMO} \|\nabla : (u^\eps - v ) \otimes (u^\eps -v)\|_{\mathcal{H}^1} ,$$ where $c'$ is a universal constant.
Then the div-curl theorem (Theorem \ref{divcurl}) applied to the $\mathcal{H}^1$ norm in the right hand side yields
\begin{eqnarray*}
 \int_{\mathbb{R}^2} v. \nabla : (u^\eps - v ) \otimes (u^\eps -v) \,dx & \leq & c \|v\|_{BMO} \|u^\eps - v\|_{L^2} \|\nabla (u^\eps -v)\|_{L^2} \\
& \leq & \frac{c^2}{2} \|v\|_{BMO}^2 \|u^\eps - v\|_{L^2}^2 + \frac{1}{2} \|\nabla (u^\eps -v)\|_{L^2}^2 \\
& \stackrel{\eqref{u-v}}{\leq} & 2c^2 \|v\|_{BMO}^2 E_{\eps, v} (t) + \frac{1}{2} \|\nabla (u^\eps -v)\|_{L^2}^2
\end{eqnarray*}
by a Young inequality. Therefore, we have showed \eqref{est1}. \qed \\
\\
\\
At this level, we have the following inequality
\begin{equation}
 E_{\eps, v} (T) - E_{\eps, v}(0) \leq 2c^2 \int_0^T \|v\|_{BMO}^2 E_{\eps, v} (t) \,dt + \mathcal{O} \left(\varepsilon ^{\nu} \right) +  \int_0^T \!\!\! \left( \eps \|\nabla :(u^\eps \otimes u^\eps )\|_{L^2}^2 - \frac{1}{2} \|\nabla (u^\eps -v)\|_{L^2}^2 \right) dt.
\end{equation}
Let us call $\tilde A^\eps$ the term that remains to estimate :
\be
\tilde A^\eps = \int_0^T A^\eps \,dt = \int_0^T  \left( \eps \|\nabla :(u^\eps \otimes u^\eps )\|_{L^2}^2 - \frac{1}{2} \|\nabla (u^\eps -v)\|_{L^2}^2 \right)\,dt .
\ee
\subsubsection{Estimating $\int_0^T A^\eps \,dt$}
The aim of this part is to prove the following lemma
\begin{lemme}
 Let $v$ and $u^\eps$ be as above. Then, for all positive $T$ and for some positive $\mu$, we have
\begin{equation}
 \int_0^T A^\eps \,dt= \int_0^T \!\!\!  \int_{\mathbb{R}^2} \left( \eps |\nabla :(u^\eps \otimes u^\eps )|^2 - \frac{1}{2} |\nabla (u^\eps -v)|^2 \right)\,dt\,dx = \mathcal{O} (\eps ^\mu).
\end{equation}
\end{lemme}
\proof First, let us write
\be
 \nabla : (u^\eps \otimes u^\eps )=u^\eps . \nabla u^\eps = u^\eps . \nabla ( u^\eps -v) + u^\eps . \nabla v 
\ee
then perform a Young inequality
\be
 | \nabla : (u^\eps \otimes u^\eps ) |^2 \leq 2 |u^\eps . \nabla ( u^\eps -v)|^2 + 2 |u^\eps . \nabla v |^2 .
\ee
So we have that
\begin{eqnarray*}
 A^\eps & = & \eps \int_{\mathbb{R}^2} |\nabla : (u^\eps \otimes u^\eps )|^2 - \frac{1}{2}  \|\nabla (u^\eps -v)\|_{L^2}^2 \\
& \leq & 2\eps \int_{\mathbb{R}^2} |u^\eps . \nabla ( u^\eps -v)|^2 + 2\eps \int_{\mathbb{R}^2} |u^\eps . \nabla v |^2 - \frac{1}{2} \|\nabla (u^\eps -v)\|_{L^2}^2  \\
& \leq & \left(2 \eps \|u^\eps\|_{L^\infty}^2 - \frac{1}{2} \right) \int_{\mathbb{R}^2} | \nabla ( u^\eps -v)|^2 + 2\eps  \int_{\mathbb{R}^2} |u^\eps . \nabla v |^2 \\
& \leq & \left( \frac{1}{2 C_1^2} - \frac{1}{2} \right) \int_{\mathbb{R}^2} | \nabla ( u^\eps -v)|^2 + 2\eps \int_{\mathbb{R}^2} |u^\eps . \nabla v |^2  \\
& \leq & 2\eps \int_{\mathbb{R}^2} |u^\eps . \nabla v |^2 ~ \leq ~ 2 \eps \|u^\eps\|_{L^\infty}^2 \int_{\mathbb{R}^2} | \nabla v |^2  ~=~ \mathcal{O}\left( \eps^s \right)
\end{eqnarray*}
since we have (see \eqref{interp_besov} and assumptions \eqref{hyp_th2d})
\be \|u_0^\eps \|_{L^\infty} \leq C_2 \|u_0^\eps \|_{\dot H^{\delta}}^{\delta} \|u_0^\eps \|_{\dot H^{1+\delta}}^{1-\delta} \lesssim \eps^{- \frac{\delta}{2}\delta} \times \eps^{\frac{s}{2} \delta} \times \eps^{-\frac{1+\delta}{2} (1-\delta)} \times \eps^{\frac{s}{2} (1-\delta)} = \eps^{\frac{s-1}{2}} . \ee
Then we have obtained the lemma, with $\mu =s$ :
\be
 \int_0^T A^\eps \,dt = \mathcal{O} (\eps ^s).
\ee
\qed 
\subsubsection{Conclusion}
Gathering the results in the previous three subsections, we obtain
\begin{equation}
 E_{\eps, v} (T) \leq 2c^2 \int_0^T \|v\|_{BMO}^2 E_{\eps, v} (t) \,dt + \mathcal{O} \left(\varepsilon ^{\nu} \right) + \mathcal{O} (\eps ^s) + E_{\eps,v}(0) .
\end{equation}
Assuming that
\be \label{hyp_energie_zero}
 \|u_0^\eps -v_0 \|_{L^2} = \mathcal{O}\left( \eps^{\frac{s}{2}}\right)
\ee
in addition to assumptions $(H)$ and \eqref{hyp}, we have
\be
 E_{\eps,v}(0)= \mathcal{O}\left( \eps^{s}\right)
\ee
by the triangle inequality. It follows that
\be
 E_{\eps,v}(T) \leq 2 c^2 \int_0^T \|v_0\|_{BMO}^2 E_{\eps,v}(t)\,dt + \mathcal{O}\left( \eps ^{\nu} \right) + \mathcal{O}\left( \eps ^{s} \right) .
\ee
Now $\nu = \dfrac{s(1+\delta)-\tilde \eps}{2(1+\delta -\tilde \eps)} \leq \dfrac{s}{2}$ since $s\leq 1$ so
\be
E_{\eps,v}(T) \leq 2 \int_0^T \|v_0\|_{BMO}^2 E_{\eps,v}(t)\,dt + \mathcal{O}\left( \eps ^{\nu} \right).
\ee
Gronwall lemma (and  $v \in L^2_T  BMO (\R^2)$) yields
\begin{equation}
 E_{\eps, v} (T) = \mathcal{O}\left( \eps^{\nu}\right)
\end{equation}
for all positive $T$ and, using inequality \eqref{u-v}
\be
 \int_{\R^2} |u^\eps -v |^2 \,dx \leq 4 E_{\eps, v}(t) ,
\ee
we deduce the convergence in the $L^\infty_T L^2(\R^2)$ norm of $u^\eps$ solution to $(NLW_\eps)$ towards $v$ solution to the $(NS)$ equations with initial data in $H^s (\R^2)^2$, where $0<s\leq 1$. \\
Theorem \ref{th2d} is now proved.

\section{The three-dimensional case : proof of Theorem \ref{th3d}}
In this section, we shall follow the plan of the previous one: we shall start by showing global existence for the damped wave equation $(NLW_\eps)$ then we will prove the convergence of this global solution $u^\eps$ towards the solution to Navier-Stokes problem with initial data in $H^{s+\frac{1}{2}}(\R^3)^3$ and $0<s<1$. Let us recall again that the critical space for $(NS)$ in $\R^3$ is $\dot H^{\frac{1}{2}}$.
\subsection{Global existence in $ \dot H^{\frac{3}{2}+\delta} \cap \dot H^{\frac{1}{2}+\delta} (\R^3, \R^3)$}
We apply the same fixed point method as in the two-dimensional case. We perform the same scale change \eqref{scaling} and retrieve system $(NLW)$ with initial data $(u_0,u_1) \in H^{\frac{3}{2}+\delta}(\R^3)^3 \times H^{\frac{1}{2}+\delta}(\R^3)^3$. \\
We conclude then that there exists a local solution $u$ to $(NLW)$, defined on the time interval $[0,T[$, for all positive real $T$ satisfying 
\begin{equation}
\label{temps_existence3d}
T  \leq \frac{1}{12 + 72 C \left( \|u_0 \|_{\dot H^{\frac{3}{2}+\delta}(\R^3)} + \|u_0 \|_{\dot H^{\frac{1}{2}+\delta}(\R^3)} + \|u_1 \|_{\dot H^{\frac{1}{2}+\delta}(\R^3)} \right) } .
\end{equation}
In particular, while $\|u (t) \|_{\dot H^{\frac{3}{2}+\delta}(\R^3)} + \|u (t) \|_{\dot H^{\frac{1}{2}+\delta}(\R^3)} + \| \partial_t u (t) \|_{\dot H^{\frac{1}{2}+\delta}(\R^3)} $ is bounded, we can reiterate the fixed point argument and extend the solution. Let $\tmax$ be the maximal existence time. We shall prove that $\tmax = + \infty $. Following the same reasoning by contradiction as in the previous section, we assume that $\tmax < + \infty$. This would imply  
\be
\|u (t) \|_{\dot H^{\frac{3}{2}+\delta}(\R^3)} + \|u (t) \|_{\dot H^{\frac{1}{2}+\delta}(\R^3)} + \| \partial_t u (t) \|_{\dot H^{\frac{1}{2}+\delta}(\R^3)} \longrightarrow + \infty ~,~~ t \rightarrow \tmax .
\ee
\\
Let us resume our wave equation $(NLW_\eps )$ with parameter $\eps$. The scaling \eqref{scaling} gives 
\be
\|u (t) \|_{\dot H^{\frac{3}{2}+\delta}} + \|u (t) \|_{\dot H^{\frac{1}{2}+\delta}} + \| \partial_t u (t) \|_{\dot H^{\frac{1}{2}+\delta}} = \eps^{\frac{\delta}{2}} \left( \sqrt{\eps} \|u^\eps \|_{\dot H^{\frac{3}{2}+\delta}} + \|u^\eps \|_{\dot H^{\frac{1}{2}+\delta}} + \eps \| \partial_t u^\eps \|_{\dot H^{\frac{1}{2}+\delta}} \right) .
\ee
Let us define, for positive $\delta$, the energy
\be
 E^{\frac{1}{2}+\delta}_\eps (t) = \int_{\mathbb{R}^3} \left( \frac{1}{2}\vert \Lambda ^{\frac{1}{2}+\delta} (u^\eps + \varepsilon \partial_t u^\eps ) \vert ^2 + \frac{\varepsilon^2}{2} \vert \Lambda^{\frac{1}{2}+\delta} \partial_t u^\eps \vert^2 + \varepsilon \vert \Lambda^{\frac{3}{2}+\delta } u^\eps \vert^2 \right) \, d x.
\ee
Then we have
\be
\eps^{\frac{\delta}{2}} \left( \sqrt{\eps} \|u^\eps \|_{\dot H^{\frac{3}{2}+\delta}(\R^3)} + \|u^\eps \|_{\dot H^{\frac{1}{2}+\delta}(\R^3)} + \eps \| \partial_t u^\eps \|_{\dot H^{\frac{1}{2}+\delta}(\R^3)} \right)  \leq C \eps^{\frac{\delta}{2}} \sqrt{E^\delta_\eps}.
\ee
If we prove that $\eps^{\frac{\delta}{2}} \sqrt{E^\delta_\eps}$ is bounded, we can deduce that the solution $u^\eps$ is global.
\paragraph{Globalization}
Let $\delta >0$ and consider the energy $E^{\frac{1}{2}+\delta}_\eps$ defined as above. We shall prove that there exists a positive constant $C_0 $ such that
\begin{equation}
 \label{controleenergie3d}
E^{\frac{1}{2}+\delta}_\eps (t) \leq C_0 \eps^{-\delta} 
\end{equation}
for all time $t>0$. In this purpose, we adapt the method used in the two-dimensional case to three space dimensions. To get around the difficulties, we use mainly interpolations and product estimates. \\
\\
First, let us recall that $\dot H^{\frac{3}{2}+\delta} \cap L^\infty (\R^3)$ is an algebra (see \cite{alinhac}) and that the inequality
\be \label{alg3d}
 \|fg\|_{\dot H^{\frac{3}{2}+\delta}} \leq C_1 \left( \|f\|_{L^\infty} \|g\|_{\dot H^{\frac{3}{2}+\delta}} + \|g\|_{L^\infty} \|f\|_{\dot H^{\frac{3}{2}+\delta}} \right)
\ee
holds for all $f,g \in \dot H^{\frac{3}{2}+\delta} \cap L^\infty (\R^3)$. 
Moreover, using the embedding $\dot B^{3/2}_{2,1} (\R^3) \subset L^\infty (\R^3)$ and the interpolation inequality
\begin{equation}
\| f\|_{\dot B^{3/2}_{2,1} (\R^3)} \leq C \|f\|_{\dot H^{\frac{1}{2}+\delta}}^{\delta}  \|f\|_{\dot H^{\frac{3}{2}+\delta}}^{1-\delta},
\end{equation}
we obtain
\begin{equation} \label{c2}
 \|f \|_{L^\infty} \leq C_2 \|f \|_{\dot H^{\frac{1}{2}+\delta}}^{\delta}  \|f \|_{\dot H^{\frac{3}{2}+\delta}}^{1-\delta}.
\end{equation}
\begin{rmq}
In $\R^3$, the Navier-Stokes problem with initial data $v_0$ in $H^{\frac{1}{2}}$ has a global solution if $\|v_0 \|_{\dot H^{\frac{1}{2}}}$ is small enough. So we could show that $E^{\frac{1}{2}+\delta}_\eps +N  E_\eps^{\frac{1}{2}}$ decreases and deduce the globality of the solution but this method requires a smallness assumption on the $\dot H^{\frac{1}{2}}$ norm of $u_0^\eps$ which depends on $\delta$ and we would have $$\|u_0^\eps\|_{\dot H^{\frac{1}{2}}} \rightarrow 0$$ when $\delta$ goes to $0$.
\end{rmq}
We shall prove that $E^{\frac{1}{2}+\delta}_\eps \left( 1+  E_\eps^{\frac{1}{2}}\right)^N$ decreases if $N$ is large enough and if the $\dot H^{\frac{1}{2}}$ norm of $u^\eps$ is small enough. First, we have to prove that $E_\eps^{\frac{1}{2}}$ decreases. To this purpose, let us compute the time derivative of the energy.
\begin{eqnarray}
 \ddt{E_\eps^{\frac{1}{2}}}(t) & = & -\int_{\R^3} \Lambda^{\frac{1}{2}} (u^\eps . \nabla u^\eps ) . \Lambda^{\frac{1}{2}} u^\eps + \int_{\R^3} \left( \eps | \Lambda^{\frac{1}{2}} (u^\eps . \nabla u^\eps ) |^2 - |\Lambda^{\frac{3}{2}}u^\eps |^2 \right) \nonumber \\
& \leq & \int_{\R^3} u^\eps . \nabla u^\eps . \Lambda u^\eps + \left( (6C_1)^2 \eps \|u^\eps\|_{L^\infty}^2 -1\right) \|u^\eps \|_{\dot H^{\frac{3}{2}}}^2 . \label{avant_est}
\end{eqnarray}
In order to estimate the first term on the right hand side, we state the following theorem:
\begin{thm} \label{juillet} Let $f \in H^{\frac{3}{2}} (\R^3)$. Then the following estimate holds.
\begin{equation}
 \label{rdvjuillet}
\vert \int_{\R^3} \Lambda f . (f.\nabla f) \vert \leq \| \Lambda^{\frac{1}{2}} f \|_{L^2 (\R^3)}  \| \Lambda^{\frac{3}{2}} f\|_{L^2 (\R^3)}^2 .
\end{equation}
\end{thm}
\proof 
Let $u \in H^{\frac{3}{2}} (\R^3)$. First, we use the duality between $\dot H^{\frac{1}{2}}(\R^3)$ and $\dot H^{-\frac{1}{2}}(\R^3)$.
\begin{equation}
 \vert \int_{\R^3} \Lambda f . (f.\nabla f) \vert  \leq  \| \Lambda f \|_{\dot H^{\frac{1}{2}}(\R^3)}  \| f.\nabla f \|_{\dot H^{-\frac{1}{2}}(\R^3)}.
\end{equation}
We shall estimate the $\dot H^{-\frac{1}{2}}(\R^3)$ norm of $f.\nabla f$. For this sake, let us notice that $f \in H^1(\R^3) \subset L^6 (\R^3)$ by a Sobolev embedding. Thus, by Hölder inequality, we have
\begin{equation}
 f.\nabla f \in L^{\frac{3}{2}}(\R^3)
\end{equation}
 since $\frac{1}{6}+\frac{1}{2}=\frac{2}{3}$. Moreover, we shall use the Sobolev embedding $H^{\frac{1}{2}}(\R^3) \subset L^3 (\R^3)$ on its dual form $L^{\frac{3}{2}} (\R^3) \subset H^{-\frac{1}{2}}(\R^3)$ .Finally, we obtain
\begin{eqnarray*}
 \vert \int_{\R^3} \Lambda f . (f.\nabla f) \vert  &\leq & \| \Lambda f \|_{\dot H^{\frac{1}{2}}(\R^3)} \| f.\nabla f \|_{\dot H^{-\frac{1}{2}}(\R^3)} \\
& \leq & \| \Lambda f \|_{\dot H^{\frac{1}{2}} (\R^3)} \| f.\nabla f \|_{ L^{\frac{3}{2}} (\R^3)} \\
& \leq &  \| \Lambda f \|_{\dot H^{\frac{1}{2}} (\R^3)} \| f\|_{L^6 (\R^3)}  \| \nabla f \|_{ L^2 (\R^3)} \\
& \leq &  \| f \|_{\dot H^{\frac{3}{2}} (\R^3)} \| f\|_{\dot H^1 (\R^3)}^2 .
\end{eqnarray*}
Now, we use the Gagliardo-Nirenberg inequality
\be 
\| f\|_{\dot H^1} \leq \sqrt{\| f \|_{\dot H^{\frac{3}{2}} } \| f\|_{\dot H^{\frac{1}{2}}}} \,.
\ee
Thus
\be
 \vert \int_{\R^3} \Lambda f . (f.\nabla f) \vert  \leq   \| f \|_{\dot H^{\frac{3}{2}} (\R^3)}^2  \| f\|_{\dot H^{\frac{1}{2}} (\R^3)} ,
\ee
which finishes the proof. \qed \\
\\
Using theorem \ref{juillet}, we shall prove the following lemma.
\begin{lemme} \label{decroissance_energie}
Let us define $0 \leq T \leq \tmax_\eps $ by
\be
 T = \sup \left\lbrace 0 \leq \tau \leq \tmax_\eps ~ :~ \forall ~ t \in [0, \tau[ , ~\|u^\eps  (t)\|_{L^\infty} < \dfrac{1}{7 C_1 \sqrt{\eps}} \right\rbrace .
\ee
Assume $\|u_0^\eps\|_{\dot H^{\frac{1}{2}}} <\dfrac{1}{16} $ and
\be \label{hypotheses_decroissance_energie}
\eps \|u_1^\eps \|_{\dot H^{\frac{1}{2}}} + \sqrt{\eps} \|u_0^\eps \|_{\dot H^{\frac{3}{2}}} = o(1)  
\ee
when $\eps$ goes to $0$. Then the energy $E^{\frac{1}{2}}_\eps$ decreases on $[0 , T ]$ and we have
\be \label{ineg_energie_0_3d}
\ddt{E_\eps^{\frac{1}{2}}}(t) \leq -\frac{1}{8} \|u^\eps (t) \|_{\dot H^{\frac{3}{2}}}^2
\ee
for all $t \in [0,T]$.
\end{lemme}
\proof Let $t \in [0,T]$. Since $\|u^\eps (t)\|_{L^\infty} <\dfrac{1}{7C_1 \sqrt{\eps}} $, inequality \eqref{avant_est} becomes
\be 
\ddt{E_\eps^{\frac{1}{2}}}(t) \leq \int_{\R^3} u^\eps . \nabla u^\eps . \Lambda u^\eps  - \frac{1}{4} \|u^\eps \|_{\dot H^{\frac{3}{2}}}^2 . 
\ee
Now, theorem \ref{juillet} (applied to first term on the right hand side) yields
\be
 \ddt{E_\eps^{\frac{1}{2}}}(t)  \leq  \|u^\eps (t)\|_{\dot H^{\frac{1}{2}}}  \|u^\eps (t) \|_{\dot H^{\frac{3}{2}}}^2 - \frac{1}{4} \|u^\eps (t)\|_{\dot H^{\frac{3}{2}}}^2 
 \leq  \left( \|u^\eps (t)\|_{\dot H^{\frac{1}{2}}} -\frac{1}{4} \right) \|u^\eps (t)\|_{\dot H^{\frac{3}{2}}}^2 .
\ee
In $t=0$, this inequality gives that $E_\eps^{\frac{1}{2}}$ decreases in a neighborhood of $0$. Let 
\be
\tau = \sup \left\lbrace 0 \leq \tilde \tau \leq T ~ :~  E_\eps^{\frac{1}{2}} ~\textrm{decreases on}~ [0, \tau] \right\rbrace .
\ee
So $0< \tau \leq T$. Assume that $\tau <T$. On the one hand, we have
\be \| u^\eps (\tau ) \|_{\dot H^{\frac{1}{2}}}^2 < 2 E_\eps^{\frac{1}{2}} (\tau) . \ee
On the other hand, for $\eps$ small enough, we have
\be E_\eps^{\frac{1}{2}}(0) \leq \|u_0^\eps \|_{\dot H^{\frac{1}{2}}}^2 + \frac{3 \eps^2}{2} \|u_1^\eps \|_{\dot H^{\frac{1}{2}}}^2  + \eps \|u_0^\eps \|_{\dot H^{\frac{3}{2}}}^2 < 2 \|u_0^\eps \|_{\dot H^{\frac{1}{2}}}^2
\ee
thanks to assumptions \eqref{hypotheses_decroissance_energie}.
Adding to that the energy decay on $[0,\tau]$, we obtain
\be \label{controle_norme_h12}
\| u^\eps (\tau ) \|_{\dot H^{\frac{1}{2}}}^2 < 2 E_\eps^{\frac{1}{2}} (\tau) \leq 2 E_\eps^{\frac{1}{2}}(0) \leq 4 \|u_0^\eps \|_{\dot H^{\frac{1}{2}}}^2
\ee
so 
\be 
\| u^\eps (\tau ) \|_{\dot H^{\frac{1}{2}}} < 2 \times \frac{1}{16} = \frac{1}{8}.
\ee
thus $E_\eps^{\frac{1}{2}}$ decreases in a neighborhood of $\tau$ and this is a contradiction with its definition.\\
Therefore $E_\eps^{\frac{1}{2}}$ decreases on $[0,T]$ and the inequality \eqref{ineg_energie_0_3d} is true for all $t \in [0,T]$.
\qed \\
\\
\begin{lemme}
 Assume $\|u_0^\eps\|_{\dot H^{\frac{1}{2}}} <\dfrac{1}{16} $ and
\be
 (H_1) ~~~~~
 i) ~ \eps^{\frac{1+\delta}{2}} \|u_0^\eps \|_{\dot{H}^{\frac{3}{2}+\delta}}  =o(1)  ~,~~ ~~~~
 ii)~ \eps^{\frac{\delta}{2}} \|u_0^\eps \|_{\dot{H}^{\frac{1}{2}+\delta}}  =o(1) .
\ee
when $\eps$ goes to zero. \\
Recall that $0 \leq T \leq \tmax_\eps $ is defined by
\be
 T = \sup \left\lbrace 0 \leq \tau \leq \tmax_\eps ~ :~ \forall ~ t \in [0, \tau[ , ~\|u^\eps  (t)\|_{L^\infty} < \dfrac{1}{7 C_1 \sqrt{\eps}} \right\rbrace .
\ee
Then $T>0$ and there exists a large number $N$, depending only on $\delta$, and a constant $C>1$ such that 
\be 
 E^{\frac{1}{2}+\delta}_\eps (t) \leq C^N \, E^{\frac{1}{2}+\delta}_\eps (0) 
\ee
for all $t\in [0,T)$ and $\eps$ small enough.
\end{lemme}
\proof
Let us compute the derivative of the energy.
\begin{eqnarray*}
 \frac{d}{dt} E^{\frac{1}{2}+\delta}_\eps  (t)&=& \int_{\R^3} \left(\Lambda^{\frac{1}{2}+\delta} (\eps \partial_{tt} u^\eps + \partial_t u^\eps ). \Lambda^{\frac{1}{2}+\delta}(u^\eps + \eps \partial_t u^\eps ) + \eps^2 \Lambda^{\frac{1}{2}+\delta} \partial_t u^\eps . \Lambda^{\frac{1}{2}+\delta} \partial_{tt} u^\eps + \right. \\
&   & ~~~ \left. + 2\eps \Lambda^{\frac{3}{2}+\delta} u^\eps . \Lambda^{\frac{3}{2}+\delta} \partial_t u^\eps \right)\\
& = &\int_{\R^3}\left( \Lambda^{\frac{1}{2}+\delta} (\eps \partial_{tt} u^\eps + \partial_t u^\eps - \Delta u^\eps ). \Lambda^{\frac{1}{2}+\delta}(u^\eps + 2\eps \partial_t u^\eps ) + \right. \\
&   & ~~~ \left.+ \Lambda^{\frac{5}{2}+\delta} u^\eps . \Lambda^{\frac{1}{2}+\delta} \left( u^\eps + 2\eps \partial_t u^\eps \right) - \eps \Lambda^{\frac{1}{2}+\delta} (\eps \partial_{tt} u^\eps + \partial_t u^\eps ). \Lambda^{\frac{1}{2}+\delta} \partial_t u^\eps + \right. \\
&   & ~~~ \left. + \eps^2 \Lambda^{\frac{1}{2}+\delta} \partial_t u^\eps . \Lambda^{\frac{1}{2}+\delta} \partial_{tt} u^\eps +2\eps \Lambda^{\frac{3}{2}+\delta} u^\eps . \Lambda^{\frac{3}{2}+\delta} \partial_t u^\eps \right).
\end{eqnarray*}
Now, recall that $u^\eps$ is a solution to $(NLW_\eps)$. Then
\begin{eqnarray*}
\frac{d}{dt} E^{\frac{1}{2}+\delta}_\eps  (t)& = &\int_{\R^3}\left( - \Lambda^{\frac{1}{2}+\delta} (\nabla : u^\eps \otimes u^\eps ). \Lambda^{\frac{1}{2}+\delta}(u^\eps + 2\eps \partial_t u^\eps ) - \Lambda^{\frac{3}{2}+\delta} u^\eps . \Lambda^{\frac{3}{2}+\delta} \left( u^\eps + 2\eps \partial_t u^\eps \right) \right. \\
&   & ~~~\left. - \eps \vert \Lambda^{\frac{1}{2}+\delta} \partial_t u^\eps \vert ^2 + 2\eps \Lambda^{\frac{3}{2}+\delta} u^\eps . \Lambda^{\frac{3}{2}+\delta} \partial_t u^\eps \right)\\
& = &\int_{\R^3} \left(- \Lambda^{\frac{1}{2}+\delta} (u^\eps . \nabla u^\eps ). \Lambda^{\frac{1}{2}+\delta}u^\eps -2\eps \Lambda^{\frac{1}{2}+\delta} (u^\eps . \nabla u^\eps ). \Lambda^{\frac{1}{2}+\delta} \partial_t u^\eps  - \right. \\
&   & ~~~\left. - \eps \vert \Lambda^{\frac{1}{2}+\delta} \partial_t u^\eps \vert ^2 - \vert \Lambda^{\frac{3}{2}+\delta} u^\eps \vert ^2 \right)\\
&\stackrel{\textrm{Young}}{\leq} & - \int_{\R^3} \Lambda^{\frac{1}{2}+\delta} (u^\eps . \nabla u^\eps ). \Lambda^{\frac{1}{2}+\delta}u^\eps + \int_{\R^3}  \left( \eps \vert \Lambda^{\frac{1}{2}+\delta} (u^\eps . \nabla u^\eps) \vert ^2 - \vert \Lambda^{\frac{3}{2}+\delta} u^\eps \vert ^2 \right) .\\
\end{eqnarray*}
Cauchy-Schwarz inequality followed by \eqref{alg3d} yields
\begin{eqnarray}
 \vert \int_{\R^3} \Lambda^{\frac{1}{2}+\delta} (u^\eps . \nabla u^\eps ). \Lambda^{\frac{1}{2}+\delta}u^\eps \vert & \leq & \| u^\eps . \nabla u^\eps \|_{\dot H^{\frac{1}{2}+\delta}}   \| u^\eps \|_{\dot H^{\frac{1}{2}+\delta}} \nonumber \\
& \leq & 3 \times 2C_1 \|u^\eps\|_{L^\infty} \|u^\eps\|_{\dot H^{\frac{3}{2}+\delta}} \| u^\eps \|_{\dot H^{\frac{1}{2}+\delta}}\nonumber \\
& \leq & 6 C_1 \|u^\eps \|_{L^\infty} \|u^\eps \|_{\dot H^{\frac{3}{2}+\delta}}   \| u^\eps \|_{\dot H^{\frac{1}{2}+\delta}} . \label{tmp}
\end{eqnarray}
Thus, back to \eqref{tmp}, we have
\begin{eqnarray*}
 \vert \int_{\R^3} \Lambda^{\frac{1}{2}+\delta} (u^\eps . \nabla u^\eps ). \Lambda^{\frac{1}{2}+\delta}u^\eps \vert & \leq & 6 C_1 \|u^\eps \|_{L^\infty} \|u^\eps \|_{\dot H^{\frac{3}{2}+\delta}}  \| u^\eps \|_{\dot H^{\frac{1}{2}+\delta}} \\
& \leq & 6 C_1 C_2  \| u^\eps \|_{\dot H^{\frac{1}{2}+\delta}}^{1+\delta} \|u^\eps \|_{\dot H^{\frac{3}{2}+\delta}}^{2-\delta} \\
& \leq & 6 C_1 C_2  \| u^\eps \|_{\dot H^{\frac{1}{2}+\delta}} \| u^\eps \|_{\dot H^{\frac{1}{2}+\delta}}^{\delta} \|u^\eps \|_{\dot H^{\frac{3}{2}+\delta}}^{2-\delta} \\
& \leq & 6 C_1 C_2 C_3 \| u^\eps \|_{\dot H^{\frac{1}{2}}}^{1-\delta} \left( \| u^\eps \|_{\dot H^{\frac{3}{2}}}  \| u^\eps \|_{\dot H^{\frac{1}{2}+\delta}} \right)^{\delta} \|u^\eps \|_{\dot H^{\frac{3}{2}+\delta}}^{2-\delta} ,
\end{eqnarray*}
where we have used another interpolation inequality between $\dot H^{\frac{1}{2}}$ and $\dot H^{\frac{3}{2}}$ :
\be
 \|f\|_{\dot H^{\frac{1}{2}+\delta}} \leq C_3 \|f\|_{\dot H^{\frac{1}{2}}}^{1-\delta} \|f\|_{\dot H^{\frac{3}{2}}}^\delta .
\ee
Furthermore, the inequality
\begin{eqnarray}
 \int_{\R^3}  \left( \eps \vert \Lambda^{\frac{1}{2}+\delta} (u^\eps . \nabla u^\eps) \vert ^2 - \vert \Lambda^{\frac{3}{2}+\delta} u^\eps \vert ^2 \right) & \leq & \left( (6 C_1)^2 \eps \| u^\eps \|_{L^\infty}^2 - 1 \right)  \int_{\R^3} \vert \Lambda^{\frac{3}{2}+\delta} u^\eps \vert ^2  \nonumber \\
& \leq & -\frac{1}{4} \|u^\eps\|_{\dot H^{\frac{3}{2}+\delta}}^2 \label{simplification}
\end{eqnarray}
holds if
\begin{equation}
\label{controlenormeinfinie3d}
 \|u^\eps (t)\|_{L^\infty} \leq \frac{1}{7 C_1 \sqrt{\eps}}.
\end{equation}
For $\eps$ small enough, this condition is ensured in $t=0$ thanks to assumptions $(H_1)$.
By continuity of the local solution $u^\eps $ with respect to the time variable $t$, we deduce that $T >0$. Let $t <T$. \\
A Young inequality yields
\begin{eqnarray*}
 \frac{d}{dt}E^{\frac{1}{2}+\delta}_\eps &\leq  &\delta 2^{2 \frac{2-\delta}{\delta}} (6 C_1 C_2 C_3 )^{\frac{2}{\delta}} \| u^\eps (t) \|_{\dot H^{\frac{1}{2}}}^{2 \frac{1-\delta}{\delta}}  \| u^\eps (t) \|_{\dot H^{\frac{3}{2}}}^2  E^{\frac{1}{2}+\delta}_\eps(t) \\
& \leq & \delta 2^{2 \frac{2-\delta}{\delta}} (6 C_1 C_2 C_3 )^{\frac{2}{\delta}} \left( 4\| u_0^\eps \|_{\dot H^{\frac{1}{2}}}^{2} \right)^{ \frac{1-\delta}{\delta}}  \| u^\eps (t) \|_{\dot H^{\frac{3}{2}}}^2  E^{\frac{1}{2}+\delta}_\eps(t)  \\
& < & 16 \delta \, (3 C_1 C_2 C_3)^{\frac{2}{\delta}}\| u^\eps (t) \|_{\dot H^{\frac{3}{2}}}^2  E^{\frac{1}{2}+\delta}_\eps(t)  .
\end{eqnarray*}
To alleviate the notations, let us set $K= K(\delta,  C_1, C_2, C_3 )=16 \delta \, (3 C_1 C_2 C_3 )^{\frac{2}{\delta}}$. 
\\
\\
As above, in subsection \ref{globalisation2d}, we prove that $E^{\frac{1}{2}+\delta}_\eps \left(1+ E^{\frac{1}{2}}_\eps \right)^N $ decreases on $[0,T)$ for $\eps$ small enough and $N$ large enough.\\
Let $0\leq t \leq T$ and $N\geq 0$. \\
We have
\beann
\ddt{} \left( E^{\frac{1}{2}+\delta}_\eps \left(1+ E^{\frac{1}{2}}_\eps \right)^N \right) (t) & = & N \ddt{E_\eps^{\frac{1}{2}}} \left(1 + E_\eps^{\frac{1}{2}} \right)^{N-1} E_\eps^{\frac{1}{2} +\delta} + \left(1 + E_\eps^{\frac{1}{2}} \right)^{N} \ddt{E_\eps^{\frac{1}{2} +\delta}}(t) \\
& \leq & \left[- \frac{N}{8} + K \left( 1 +  E_\eps^{\frac{1}{2}} (t)\right) \right] \| u^\eps (t) \|_{\dot H^{\frac{3}{2}}}^2 \left(1 + E_\eps^{\frac{1}{2}}(t) \right)^{N-1} E_\eps^{\frac{1}{2} +\delta}(t) \\
& \leq & \left[- \frac{N}{8} + K \left( 1 +  2 \|u_0^\eps\|_{\dot H^{\frac{1}{2}}}^2 \right) \right] \| u^\eps (t) \|_{\dot H^{\frac{3}{2}}}^2 \left(1 + E_\eps^{\frac{1}{2}}(t) \right)^{N-1} E_\eps^{\frac{1}{2} +\delta}(t).
\eeann
Inequality \eqref{controle_norme_h12} yields
\be
\ddt{} \left( E^{\frac{1}{2}+\delta}_\eps \left(1+ E^{\frac{1}{2}}_\eps \right)^N \right) (t) \leq \left(- \frac{N}{8} + \frac{129}{128} K  \right) \| u^\eps (t) \|_{\dot H^{\frac{3}{2}}}^2 \left(1 + E_\eps^{\frac{1}{2}}(t) \right)^{N-1} E_\eps^{\frac{1}{2} +\delta}(t).
\ee
Now, taking $N \geq 8 K \times \frac{129}{128} = 129 \delta \, (3 C_1 C_2 C_3)^{\frac{2}{\delta}}$, we prove that $E^{\frac{1}{2}+\delta}_\eps \left(1+ E^{\frac{1}{2}}_\eps \right)^N$ decays on $[0,T]$. So
\be
E^{\frac{1}{2}+\delta}_\eps (t) \leq E^{\frac{1}{2}+\delta}_\eps (0) \left(1+ E^{\frac{1}{2}}_\eps (0) \right)^N \leq E^{\frac{1}{2}+\delta}_\eps (0) \left( \frac{129}{128} \right)^N  
\ee
for $t\in[0,T]$. \qed \\
\begin{lemme} 
Assume the limit 
\be \label{hyp3d}
 \eps^{1+\frac{\delta}{2}} \|u^\eps_1\|_{\dot H^{\frac{1}{2}+\delta}} \longrightarrow 0~, ~~\eps \rightarrow 0
\ee
in addition to the assumptions $(H_1)$ and \eqref{hypotheses_decroissance_energie}. Then we have
\begin{equation}
 E^{\frac{1}{2}+\delta}_\eps (0) =o \left( \eps^{- \delta} \right).
\end{equation}
Moreover, for all $t \in [0,T[$ and $\eps$ small enough, the inequality
\begin{equation} 
\|u^\eps \|_{L^\infty} ~\leq ~ \frac{1}{14 C_1 \sqrt{\eps}} 
\end{equation}
holds.
\end{lemme}
We skip the proof here as it is analogous to the one of Lemma \ref{lemme2}. As above, this lemma implies the inequality \eqref{controleenergie3d} and the globality of the solution $u^\eps$.
\newcommand{\E}{\mathcal{E}}
\subsection{Convergence towards a solution to $(NS)$ problem}
The aim of this subsection is to prove an error estimate in the $L^\infty_T \dot H^{\frac{1}{2}}(\R^3)^3$ norm. More precisely, let $v_0 \in H^{s+\frac{1}{2}}(\R^3)^3$, $0<s<1$, be the initial data for the Navier-Stokes equations \eqref{ns} and $(u_0^\eps , u_1^\eps) \in H^{\frac{3}{2}+\delta} \times H^{\frac{1}{2}+\delta} (\R^3)^3$, $0<\delta<1$, be the initial data for the damped wave equation \eqref{nlw}. We will prove that, if $\|u_0^\eps -v_0 \|_{\dot H^{\frac{1}{2}}}^2 = \mathcal{O}(\eps^{s})$, then for all positive $T$, we prove that
\be
\sup_{t\in (0,T)} \|(u^\eps - v)(t)\|_{\dot H^{\frac{1}{2}}}^2 = \mathcal{O}\left(\eps^{\left(\frac{s}{2}\right)^-}\right).
\ee

In the following, we adapt the method in section \ref{convergence2d}. So let us consider again the Dafermos modulated energy defined by
\begin{equation} \label{dafermos3d}
 \E_{\eps, v} (t) = \int_{\R^3} \left( \frac{1}{2} |\Lambda^{\frac{1}{2}} \left(u^\eps -v(t,x) +\eps \partial_t u^\eps \right) |^2 + \frac{\eps^2}{2} |\Lambda^{\frac{1}{2}}  \partial_t u^\eps |^2 + \eps |\Lambda^{\frac{1}{2}} \nabla u^\eps |^2 \right) \,dx
\end{equation}
for all divergence-free vector $v$. As above, this energy satisfies the inequality
\begin{equation} \label{u-v3d}
 \|u^\eps (t) -v(t) \|_{\dot H^{\frac{1}{2}}}^2 \leq 4 \E_{\eps, v} (t).
\end{equation}
We shall show that the assumptions \eqref{hyp_th3d} in Theorem \ref{th3d} and a Gronwall inequality imply that
\be \forall~  t>0~,~~ \E_{\eps,v}(t) = \mathcal{O}\left( \eps^{\left(\frac{s}{2}\right)^-} \right) . \ee
Let us compute the derivative of $\E_{\eps,v}$.
\newcommand{\leray}{\mathbb{P}}
\begin{lemme}
 The Dafermos modulated energy $\E_{\eps ,v}$ satisfies the identity
\begin{eqnarray}
 \frac{d}{dt} \E_{\eps,v}(t) & = &  \int_{\R^3} \Lambda (u^\eps -v) .  \leray \left(\nabla : u^\eps \otimes u^\eps - \nabla : v \otimes v \right) \,dx -\eps \int_{\R^3} |\Lambda^{\frac{1}{2}} \left(\partial_t u^\eps + \leray \nabla : u^\eps \otimes u^\eps \right)|^2 \,dx \nonumber \\
& & - \eps \int_{\R^3}  \Lambda^{\frac{1}{2}} \partial_t v . \Lambda^{\frac{1}{2}} \partial_t u^\eps \,dx + \int_{\R^3} \left( \eps |\Lambda^{\frac{1}{2}} \leray \nabla : (u^\eps \otimes u^\eps)|^2 - |\Lambda^{\frac{3}{2}} (u^\eps -v)|^2 \right)\,dx \nonumber \\
& & - \int_{\R^3} \Lambda (\partial_t v + \leray v.\nabla v -\Delta v).(v-u^\eps)\,dx . \label{soisnul3d}
\end{eqnarray}
\end{lemme}
As above in section \ref{convergence2d}, we take $v$ solution to $\partial_t v + \leray v.\nabla v -\Delta v =0$ in the following, so that the last term in identity \eqref{soisnul3d} vanishes. Moreover, the second term in \eqref{soisnul3d} being negative, we have
\begin{eqnarray}
 \frac{d}{dt} \E_{\eps,v}(t) & \leq & - \eps \int_{\R^3} \Lambda^{\frac{1}{2}} \partial_t v . \Lambda^{\frac{1}{2}} \partial_t u^\eps \,dx + \int_{\R^3} \Lambda (u^\eps -v). \leray \nabla : (u^\eps \otimes u^\eps - v \otimes v ) \,dx    \nonumber \\
& &  + \int \left( \eps |\Lambda^{\frac{1}{2}} \leray \nabla : (u^\eps \otimes u^\eps)|^2 - |\Lambda^{\frac{3}{2}} (u^\eps -v)|^2 \right). \label{derivee_dafermos3d}
\end{eqnarray}
Concerning the first term, the argument in subsection \ref{partie1convergence2d} applies and immediately yields the estimate
\be
\varepsilon \int_0^T\!\!\! \int_{\mathbb{R}^3} \Lambda^{\frac{1}{2}} \partial_t u^\eps . \Lambda^{\frac{1}{2}} \partial_t v \,dt\,dx = \mathcal{O}\left( \eps^{\left(\frac{s}{2}\right)^-} \right).
\ee
We shall estimate the two remaining terms on the right hand side in the following two subsections.
\subsubsection{Estimating $ \int_0^T\!\!\! \int_{\mathbb{R}^3} \Lambda (u^\eps-v). \nabla : (u^\eps \otimes u^\eps -v \otimes v ) \,dt\,dx $}
First, let us recall the identity
\be
\nabla : (u^\eps \otimes u^\eps - v \otimes v) = (v-u^\eps).\nabla v + u^\eps . \nabla(v-u^\eps) .
\ee
Using a variant of theorem \ref{juillet} above, we obtain
\begin{eqnarray*}
 -\int_{\mathbb{R}^3} \Lambda (u^\eps-v). \nabla : (u^\eps \otimes u^\eps -v \otimes v ) \,dx & = &  \int_{\mathbb{R}^3} \Lambda (u^\eps-v). \left( (v-u^\eps). \nabla v + u^\eps . \nabla (v-u^\eps) \right) \,dx \\
& \lesssim & \|u^\eps -v\|_{\dot H^{\frac{3}{2}}} \left( \|u^\eps - v\|_{L^6} \|v\|_{\dot H^1} + \|u^\eps\|_{L^6} \|u^\eps - v\|_{\dot H^1} \right) \\
& \lesssim & \|u^\eps -v\|_{\dot H^{\frac{3}{2}}} \|u^\eps - v\|_{\dot H^1} \left( \|v\|_{\dot H^1} + \|u^\eps\|_{\dot H^1} \right).
\end{eqnarray*}
Now, a Galgliardo-Nirenberg inequality followed by a Young inequality yields
\begin{eqnarray*}
 \left|\int_{\mathbb{R}^3} \Lambda (u^\eps-v). \nabla : (u^\eps \otimes u^\eps -v \otimes v ) \,dx \right| & \leq & C \|u^\eps -v\|_{\dot H^{\frac{3}{2}}}^{\frac{3}{2}} \|u^\eps - v\|_{\dot H^\frac{1}{2}}^{\frac{1}{2}} \left( \|v\|_{\dot H^1} + \|u^\eps\|_{\dot H^1} \right)\\
& \leq & \frac{C^4}{4} \|u^\eps - v\|_{\dot H^\frac{1}{2}}^2 \left( \|v\|_{\dot H^1} + \|u^\eps\|_{\dot H^1} \right)^4  + \frac{3}{4} \|u^\eps -v\|_{\dot H^\frac{3}{2}}^2 \\
& \leq & 2 C^4 \|u^\eps - v\|_{\dot H^\frac{1}{2}}^2 \left( \|v\|_{\dot H^1}^4 + \|u^\eps\|_{\dot H^1}^4 \right)  + \frac{3}{4} \|u^\eps -v\|_{\dot H^\frac{3}{2}}^2.
\end{eqnarray*}
Finally, using inequality \eqref{u-v3d}, we obtain
$$
\left|\int_{\mathbb{R}^3} \Lambda (u^\eps-v). \nabla : (u^\eps \otimes u^\eps -v \otimes v ) \,dx \right| \leq 8C^4 \left( \|v\|_{\dot H^1}^4 + \|u^\eps\|_{\dot H^1}^4\right) \E_{\eps, v} (t) + \frac{3}{4} \|u^\eps -v\|_{\dot H^\frac{3}{2}}^2.
$$
Now, let $T>0$ and recall that $v$ is a solution to the Navier-Stokes equations, with initial data $v_0 \in H^{s+\frac{1}{2}} (\mathbb{R}^3)^3$ and $0<s<1$. So $v \in L^2_T  H^{s+\frac{3}{2}}\cap L^\infty_T H^{s+\frac{1}{2}} (\mathbb{R}^3)^3$. By interpolation, we have also $v \in L^4_T  H^{s+1}(\mathbb{R}^3)^3$.\\
Furthermore, notice that inequality \eqref{controle_norme_h12} implies that  $u^\eps \in L^\infty_T \dot H^{\frac{1}{2}}$ uniformly in $\eps$ and, using the energy decay proven in lemma \ref{decroissance_energie}, one can easily show that $u^\eps \in L^2_T \dot H^{\frac{3}{2}}$ uniformly in $\eps$. Thus, by interpolation, we obtain that 
$$ 
u^\eps \in L^4_T \dot H^{1}  (\mathbb{R}^3)^3
$$ uniformly in $\eps$. Therefore
\begin{eqnarray}
\int_0^T\!\!\! \int_{\mathbb{R}^3}  \Lambda (u^\eps-v). \nabla : (u^\eps \otimes u^\eps -v \otimes v ) \,dt\,dx & \leq & 8C^4 \int_0^T \|v (t) \|_{H^{s+1}}^4 \E_{\eps, v} (t) \,dt \nonumber \\
& & + 8C^4 \int_0^T \|u^\eps (t) \|_{\dot H^{1}}^4 \E_{\eps, v} (t) \,dt \nonumber \\
& & + \frac{3}{4} \int_0^T \|\nabla (u^\eps -v)\|_{\dot H^\frac{1}{2}}^2 . \label{est13d}
\end{eqnarray}
At this level, we have the following estimate:
\begin{eqnarray}
\E_{\eps, v} (T) - \E_{\eps,v}(0) &\leq & 8C^4 \int_0^T \left( \|v\|_{H^{s+1}}^4 + \|u^\eps \|_{\dot H^1}^4 \right) \E_{\eps, v} (t) \,dt + \mathcal{O} \left(\eps^{\left(\frac{s}{2}\right)^-} \right) \nonumber \\
& & + \int_0^T \!\!\!  \left( \eps \|\nabla: u^\eps \otimes u^\eps\|_{\dot H^{\frac{1}{2}}}^2 - \frac{1}{4}  \|\nabla (u^\eps -v) \|_{ \dot H^{\frac{1}{2}} }^2 \right) (t)\, dt.
\end{eqnarray}
Let us set
\be \label{Aeps}
 \int_0^T A^\eps (t) \,dt := \int_0^T \!\!\!  \left( \eps \|\nabla: u^\eps \otimes u^\eps\|_{\dot H^{\frac{1}{2}}}^2 - \frac{1}{4} \|\nabla (u^\eps -v) \|_{ \dot H^{\frac{1}{2}} }^2 \right) (t)\, dt.
\ee
It remains to estimate the last term $ \int_0^T A^\eps (t) \, dt$.
\subsubsection{Estimating $\int_0^T A^{\eps} (t) \,dt$}
First, let us recall that $L^\infty \cap \dot H^{\frac{3}{2}}(\R^3)$ is an algebra. So 
\be
 \| \nabla :(u^\eps \otimes u^\eps )\|_{\dot H^{\frac{1}{2}}} \lesssim \| u^\eps \otimes u^\eps \|_{\dot H^{\frac{3}{2}}} \leq C \|u^\eps \|_{L^\infty} \|u^\eps \|_{\dot H^{\frac{3}{2}}}.
\ee
Therefore we have
 \beann
 A^\eps (t) & \leq & C^2 \eps \|u^\eps \|_{L^\infty}^2 \|u^\eps \|_{\dot H^{\frac{3}{2}}}^2 - \frac{1}{4} \|u^\eps -v \|_{ \dot H^{\frac{3}{2}} }^2 \\
& \leq & \left( 2C^2 \eps \|u^\eps \|_{L^\infty}^2 - \frac{1}{4} \right) \|u^\eps -v \|_{ \dot H^{\frac{3}{2}} }^2 + 2C^2 \eps \|u^\eps \|_{L^\infty}^2 \|v \|_{\dot H^{\frac{3}{2}}}^2
\eeann
by a Young inequality. 
Now, using that $\eps \|u^\eps \|_{L^\infty\left((0,T)\times \R^3\right)}^2 = \mathcal{O}\left(\eps^s \right)$ and that $v\in L^2_T \dot H^{\frac{3}{2}}$, we obtain
\be
 \int_0^T A^\eps (t) \,dt \leq 2C \eps \|u^\eps \|_{L^\infty}^2 \| v\|_{L^2_T \dot H^{\frac{3}{2}}}^2 =\mathcal{O}\left(\eps^s \right)
\ee
if $\eps$ is small enough.

\subsubsection{Conclusion}
Finally, we have 
\begin{equation}
\E_{\eps, v} (T) \leq 8C^4 \int_0^T \left( \|v\|_{H^{s+1}}^4 + \|u^\eps \|_{\dot H^1}^4 \right) \E_{\eps, v} (t) \,dt + \mathcal{O} \left(\eps^{\left(\frac{s}{2}\right)^-} \right) + \E_{\eps,v}(0).
\end{equation}
Assuming that 
\be
\|u_0^\eps -v_0 \|_{ \dot H^{\frac{1}{2}} } + \eps \|u_1^\eps \|_{ \dot H^{\frac{1}{2}} } + \sqrt{\eps} \|u_0^\eps \|_{ \dot H^{\frac{3}{2}} } = \mathcal{O}\left( \eps^{\frac{s}{2}}\right) ,
\ee
we have immediately
\be
\E_{\eps ,v}(0) = \mathcal{O}\left(\eps^s \right) .
\ee
Let us recall that $v \in L^4_T   H^{s+1} (\R^3)^3$ and $u^\eps \in L^4_T \dot H^1$. Gronwall lemma implies that 
\begin{equation}
 \E_{\eps, v} (T) = \mathcal{O} \left(\eps ^{ \left(\frac{s}{2}\right)^- }\right)  . 
\end{equation}
Now, \eqref{u-v3d} yields the announced error estimate and the proof of Theorem \ref{th3d} is complete.

\section{Aknowledgements}
I would like to thank my supervisors Valeria Banica and Pierre-Gilles Lemarié-Rieusset from the Laboratoire Analyse et Probabilités at the Université d'Evry-Val d'Essonne for their valuable discussions, help and advice.
\bibliographystyle{unsrt}

\end{document}